\documentclass[12pt]{article}
\usepackage[pdftex]{graphicx}
\usepackage{amsfonts}
\usepackage{mathrsfs}
\newtheorem{theo}{Theorem}

\newtheorem{prop}{Proposition}
\newtheorem{lemm}{Lemma}

\newtheorem{remk}{Remark}

    \newcommand{\eqco}{\setcounter{equation}{0}}
    \newcommand{\thco}{\setcounter{theo}{0}}
    \newcommand{\prco}{\setcounter{prop}{0}}
    \newcommand{\laco}{\setcounter{lemm}{0}}
    \newcommand{\coco}{\setcounter{coro}{0}}
    \newcommand{\cjco}{\setcounter{conj}{0}}
    
    \newcommand{\deco}{\setcounter{defn}{0}}
    
    \newcommand{\allco}{\eqco  \thco \prco \laco \coco \cjco \deco}
\newcommand{\bea}{\begin{eqnarray}}
\newcommand{\eea}{\end{eqnarray}}
\newcommand{\beaa}{\begin{eqnarray*}}
\newcommand{\eeaa}{\end{eqnarray*}}
\newcommand{\Z}{{\mathbb Z}}

\newcommand{\N}{{\mathbb N}}

\newcommand{\lbl}{\label}
\newcommand{\bean}{\begin{eqnarray*}}
\newcommand{\eean}{\end{eqnarray*}}
\renewcommand{\epsilon}{\varepsilon}
\newcommand{\eps}{\varepsilon}

\newcommand{\ty}{\tilde{y}}

\title{{Percolation of even sites for random sequential adsorption}}
\author{ Mathew D. Penrose\footnote{ Department of Mathematical Sciences,
 University of Bath, Bath BA2 7AY, United Kingdom,
Email: m.d.penrose@bath.ac.uk}
~ and Tom Rosoman\footnote{ Department of Mathematical Sciences,
 University of Bath BA2 7AY, United Kingdom,
Email: ter20@bath.ac.uk}  \footnote{ Supported by an EPSRC studentship}
\\
{\em University of Bath}
 }

\date{\today}

\begin{document}

\maketitle

\begin{abstract}\addcontentsline{toc}{section}{Abstract}
Consider random sequential adsorption on a red/blue chequerboard 
lattice with arrivals at rate $1$ on the red squares and rate $\lambda$ on
 the blue squares. We prove that the
 critical value of $\lambda$, above which we get
 an infinite blue component, is finite and strictly greater than $1$.

\end{abstract}



%\section{Introduction and statement of result}
\section{Introduction}
\lbl{secintro}
Random sequential adsorption (RSA) is a term used for a
family of probabilisitic models for irreversible
particle deposition.
Particles arrive at random locations onto a surface
which is typically taken to be two-dimensional and
 initially empty, and each
particle, once accepted, blocks nearby locations from
becoming occupied, thereby causing any subsequent particles
arriving nearby to be rejected. Both lattice and continuum versions
of  RSA have been studied extensively in the literature.
They are of considerable interest in the physical sciences,
for example with regard to the coating of a surface by
some  adsorbed substance \cite{ev93,Privman}.

In the present paper we consider the two-dimensional
lattice version of RSA, whereby the surface
is represented by
 $\Z^2$ endowed with the usual nearest-neighbour
graph structure. The arrival times $t_x, x \in \Z^2,$
are taken to be independent and exponentially
distributed. Initially all sites are vacant,
but if a particle arrives at $x$ at time
$t_x$, then site $x$ becomes occupied at that instant
 unless one of the neighbouring sites was previously
occupied. That is, when a particle becomes occupied
it causes all of its neighbours to be blocked. Ultimately,
every site is either occupied or blocked. Provided
there is a uniform bound on the arrival rates,
the model is well defined even on the infinite
lattice $\Z^2$: see e.g. \cite{Pen02} or \cite{PS}. 

The ultimate  configuration is
called the {\em jammed state} and is the focus
of our attention here.  In the jammed state,
the occupied lattice sites comprise 
a maximal stable set (a stable set
is a subset of vertices in the graph such that no two
vertices in that set are adjacent). The remaining sites
are blocked. 

Since the $\Z^2$ lattice is bipartite,
the set of occupied sites is naturally partitioned
into two phases, the even and odd occupied sites,
where a site is denoted even/odd according to its graph
distance from the origin. In fact, we can partition
the whole of $\Z^2$ into two phases,
one phase consisting of even occupied sites and odd blocked sites
(the {\em even phase}),
and the other consisting of
odd occupied sites and even  blocked sites (the {\em odd phase}).
Since we are in the jammed state,  all sites lie in one phase or the other.
 
We are interested in the {\em percolation}  properties of the 
even phase.  That is, we consider the question of whether
the subgraph of $\Z^2$ induced by the even phase contains
an infinite component.
  Physically, such questions
could be of interest with regard to, for example, electrical
or thermal conductivity through adsorbed particles on a surface.
Percolation properties of particle configurations
generated by RSA type processes have been studied
in the physical sciences literature; see for example
 Section VI of \cite{ev93}, \cite{Loscar} and \cite{Rampf}, and
 references therein. 

The sites in the even phase form a dependent site percolation
process on $\Z^2$.
%By analogy with known results for independent  percolation,
%one might expect that if the arrival rates at all sites 
%are the same, then the even phase will not percolate (and 
%neither will the odd phase). 
A basic result of this paper is that
if the arrival rates at all sites 
are the same, then the even phase will not percolate (and 
neither will the odd phase).  Therefore
%As we shall show, this turns out to be the case,
%and therefore,
 the odd and even phases decompose
into finite connected islands (cf. the diagrams
on page 1309 of \cite{ev93}).

One can tune the model by biasing the arrival rates
in favour of the even sites, and this is what we do, with
a single parameter $\lambda$  representing the amount of bias
(this version of the model was suggested to us by Martin Zerner).
One might expect the even phase to percolate given
a sufficiently high level of bias.
We shall show that there is a non-trivial phase transition 
in the parameter $\lambda$. In particular, there is a non-zero
level of bias at which the even phase still does not percolate.
This improves on  
the aforementioned basic result, and
is our main result.
% and is stronger than the
% improving the weak inequality
%for the critical value of $\lambda$ implied by our basic result

We briefly discuss the degree of surprise
in the basic result of non-percolation when all arrival rates
are the same. It is known that independent site percolation with parameter
$p=1/2$ on the usual square lattice does not percolate, and
the density of the even phase in RSA is 1/2, suggesting by
analogy that our  dependent site percolation process
would  not percolate.  On the other hand, if one turns
the lattice through 45 degrees and 
considers only the even sites,
 these form  a square grid
with diagonal connections regarding occupied sites,
 since for any two occupied even sites two steps apart, the intervening
odd site must be vacant. Site percolation on  the square lattice with diagonals
is strictly supercritical at $p=1/2$  (it is dual to
site percolation on the usual square lattice
which is strictly subcritical), so from this one might 
expect the even phase  of a RSA-type hard-core process
(i.e., one which generates a random stable subset of $\Z^2$)
 with density sufficiently
close to $1/2$  to percolate.

It seems that the second of these two analogies is misleading here.
Indeed, it seems likely that RSA can be modified 
to provide a hard-core process with
 density of occupied sites arbitrarily 
close to one-half, without affecting the 
 basic non-percolation result. To see this, consider
a variant where, initially, large square blocks 
%$(K_1 \times K_1)$
of sites arrive sequentially at random locations.
When a square block arrives, suppose all sites in the block with
the same parity as its lower left corner become occupied,
unless one or more of them is already blocked, in which
case the entire incoming square block is rejected.
At the end of this process there remain some holes,
 but these can be filled in by having
a subsequent arrivals process of smaller square blocks.

On the other hand, we think it is also likely that there
exist stationary ergodic hard-core processes on the
sites of the square lattice for which the even phase does percolate. 
Indeed, consider a  stationary curve along
the lines of  the one in the proof of Proposition 5 of
 Holroyd and Liggett \cite{HoLig},
and put the odd phase on one side of this and
the even phase on the other side.  

Van den Berg \cite{Berg11} considers another form dependent percolation,
 with biological motivation.  That paper
is concerned with sharp transitions for percolation on
 a the random field associated
with the contact process, whereas in the present instance we
 are concerned with inequalities  of critical points
for a random field generated by random sequential adsorption.
It is noteworthy, however, that in both cases the
methods of Bollob\'as and Riordan \cite{BR06} play a key role.

\section{Statement of result}

We now describe the model in more detail.
 We 
have a chequerboard pattern of sites, where
 each even site is adjacent to four odd sites and vice versa.
 Assume we have a Poisson arrival process at each site
 with rate $1$ on the odd sites and rate $\lambda$ on the even sites
(these arrival processes are taken to be mutually independent).
 We start off with all sites empty. 
Let $t_x$ be the time of the first arrival in the Poisson process at $x$.
If none of the four neighbours of $x$ are occupied at time $t_x$ we declare
$x$ to be occupied from then on. 
If any of its neighbours becomes occupied then site $x$ becomes blocked
at that instant, and remains so from then on.
 In this way every site will eventually end up being occupied or blocked (see
\cite{Pen02} or \cite{PS}.)

If an even site is occupied we declare it to be black and if it is blocked we declare it to be white. If an odd site is occupied we declare it to be white and if it is blocked we declare it to be black. The black sites form the even
 phase mentioned earlier.
 We form a graph of black vertices with edges between any two black
 vertices that are adjacent in the square lattice.
By the ergodic property of any family of independent identically
distributed variables indexed by $\Z^2$, the probability that
there is an infinite black component is either zero or one. Moreoever,
by a standard coupling argument, this probability is monotonic
nondecreasing in $\lambda$. Therefore,
 there is a critical value $\lambda_c \in [0,\infty]$, such
that for $\lambda > \lambda_c$ there will
 almost surely be an infinite black component and for
$\lambda < \lambda_c$ there
 will almost surely not be an infinite black component.
Our main result provides some non-trivial bounds on this critical value.

\begin{figure}[htbp]
\begin{center}
\scalebox{.6}{\includegraphics[angle = 270]{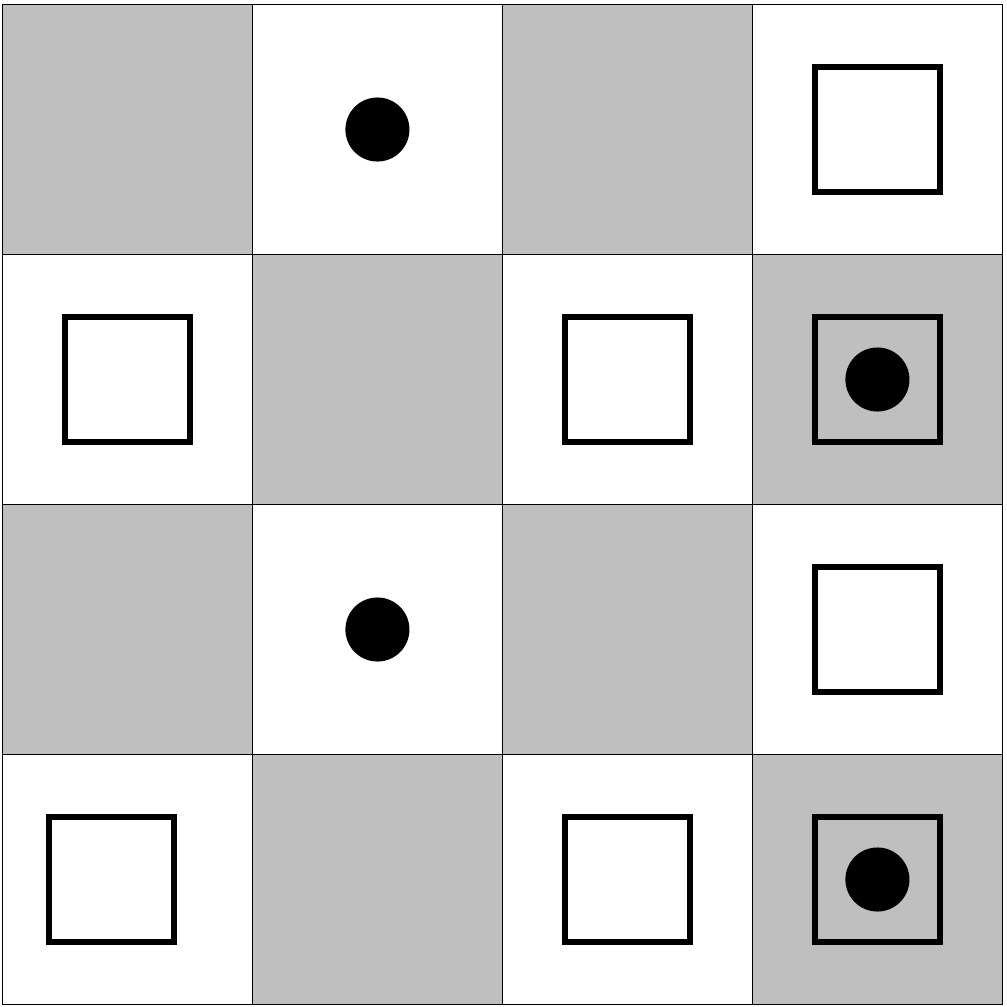}}
\end{center}
\caption{Example: The shaded squares are even sites and the unshaded
squares are odd sites. The squares with a circle in are occupied sites,
 the squares with an inscribed square in are black sites and
 the squares without an inscribed square in are white sites}
\label{fig:enh1}
\end{figure}

%\section{Our Result}
\begin{theo}
\label{thm1}
It is the case that
$1 < \lambda_c < 10$.
\end{theo}
{\em Proof of  $\lambda_c < 10$.}
This upper bound is simple to prove, and we deal with it at once. 
 We start by colouring all sites yellow that have 
 even coordinates adding up to a multiple of $4$,
such as $(0,0), (2,2), (0,4)$ and so on.
 We define a square lattice of yellow sites by saying two yellow
 sites are adjacent if they
are $2\sqrt{2}$ apart. We then consider site percolation on this lattice with each site being occupied independently with probability $\frac{\lambda}{4 + \lambda}$. This corresponds to the probability that an arrival at a yellow site happens before an arrival at any of its neighbours in the original lattice. If two adjacent yellow sites are occupied in the new lattice then they are occupied in the original lattice and also the even site midway between them will also be occupied. Therefore if there is an infinite component in the new lattice there is also one in the original lattice, so we have the following inequality:
\[
 \frac{\lambda_c}{4 + \lambda_c} \leq p_s
\]
where $p_s$ is the critical site probability on the square lattice, which is
known to be less than $0.7$ (Wierman \cite{Wier95}). Rearranging gives that
\[
\lambda_c \leq \frac{4 p_s}{1 - p_s} < \frac{28}{3} < 10
\]
so we have proved the upper bound.
\hfill{$\Box$} \\

In the remaining sections, we shall prove the lower bound
$\lambda_c >1$. Although the result
is perhaps to be expected by analogy with known (though non-trivial)
results for Bernoulli (i.e., independent)
site percolation, we are not aware
of any such results in a dependent site percolation
setting such as we consider here.
% and our proof is somewhat technical.
% Our approach requires some fairly sophisticated
%techniques of modern percolation theory.
 By use of the weak RSW-type lemma
established by Bollob\'as and Riordan \cite{BR06} for
 percolative systems enjoying weak dependence, 
we shall rather quickly establish the weak
 version of the inequality, namely
$\lambda_c \geq 1$ (see Remark \ref{weakineq}).
To make this inequality strict we use 
 the technique of enhancement. While this technique is well known,
in the present setting  its application is quite intricate,
requiring a whole sequence of notions of pivotal vertex 
(see Sections \ref{secenhance} and \ref{secpivcomp}).

\section{Duality}
\label{secdual}
\eqco
Define the dual lattice $\Lambda^*$ to be the square lattice
 $\Lambda$ with the diagonals added so that two sites are adjacent
 if their centres are at distance $1$ or $\sqrt{2}$ from each other.
 On any square set of sites we have exactly one of the following two
 events, either a black horizontal crossing in $\Lambda$ or a
 white vertical crossing in $\Lambda^*$. 

Define $f_{\lambda}(\rho, s)$ to be the probability that there is
a horizontal black crossing of the rectangle
$[1,2 \lfloor  \frac{\rho s}{2} \rfloor ] \times
[1,2 \lfloor \frac{s}{2} \rfloor ] $ (an approximately
$\rho s \times s$ lattice rectangle with even side lengths).
 Define $f^*_{\lambda}(\rho, s)$ to be the probability that
there is a horizontal black crossing of 
this  rectangle when
% a $\rho s$ by $s$ rectangle when
%
% where we take the largest rectangle with all sides of
% even length that fits inside it and where 
we allow diagonal edges as well. 

In subsequent sections, we shall prove the following key result.
\begin{prop}
\lbl{keyprop}
There exists $\mu < 1$ such that 
\begin{equation}
\liminf_{s \rightarrow \infty}f^*_{\mu}(1,s) > 0.
\label{eqstar}
\end{equation}
\end{prop}
In the remainder of the present section, we show how to complete
the proof of Theorem \ref{thm1},  given Proposition \ref{keyprop}.
The argument uses two further results, which we give now.

We say site $x \in \Z^2$ {\em affects} site $y \in \Z^2$ if there exists
a self-avoiding path in $\Z^2$  starting at $x$ and ending at $y$,
with arrival times occurring in increasing order  
along this path. If $x$ does not affect $y$, then any change
to $t_x$ (with other arrival times unchanged)
will not cause any change to the occupied/blocked status of site
$y$.
Similarly to arguments in \cite{Pen02}, we have
the following simple lemma.
\begin{lemm}
\lbl{affectlem}
Let $x \in  \Z^2$.
 The probability that site $x$ is affected
from distance greater than $r$ tends to zero as $r \to \infty$.
 Likewise, the probability that site $x$ affects 
some site at distance
greater than $r$ from $x$ tends to zero as $r \to \infty$.
\end{lemm}
{\bf Proof.} 
For any self-avoiding path of length $r$, taking alternate sites along the path
one has at least $\lfloor r/2 \rfloor$ 
independent identically  distributed arrival times, so the probability
they occur in increasing order is  at most
 $1/\lfloor r/2 \rfloor!$. 
Therefore the probability that $x$
% a  particular vertex in $B(3n)$ is
is affected from distance greater than $r$ 
is at most
%$\frac{4*3^m}{m!}$
$4(3^r)/\lfloor r/2 \rfloor!$, 
 which tends to zero as $r \to \infty$.
The proof of the second part is  similar.
\hfill{$\Box$} \\

We also use the following much deeper lemma, which is a weak version of the
 RSW lemma for dependent percolation.
\begin{prop}
\label{RSWlem}
Let $\lambda > 0$ and $\rho > 1$ be fixed. If $\liminf_{s \rightarrow \infty}f^*_{\lambda}(1,s) > 0$ then $\limsup_{s \rightarrow \infty} f^*_{\lambda}(\rho,s) > 0$. 
\end{prop}
A result along these lines is given by
  Bollob\'as and Riordan (Theorem 4.1 of  \cite{BR06}).
 The result in \cite{BR06} is for Voronoi percolation but the
 proof can be transferred to our model,
as we now discuss.

Much of the proof in \cite{BR06} 
 relies only on the Harris-FKG inequality, which
 holds in the present  model as well (see 
Penrose and Sudbury \cite{PS}).
 These arguments in \cite{BR06} carry over easily, making sure
 that rectangles with even integer sides are chosen as the RSA model
 is on a discrete lattice not a continuum.

In the first part of the proof in \cite{BR06}, an event
 $E_{dense}$ is considered, and we need a 
 different version of this event here.
Given an integer $s$ and constant $\rho > 1$ let $R_s$ be an $s$ by $ \lfloor \rho s \rfloor  $ rectangle. Given a rectangle $R$ with integer sides $a$ and $b$ let $R[r]$ be the rectangle with sides $a + 2r$ and $b +2r$ centred on $R$, so the edges of $R[r]$ are at distance $r$ from the edges of $R$. Let $E_{dense}(R_s)$ be the event that no site in $R_s$
is affected by any site outside  outside $R_s[r-1]$,
 where we take $r$ to be $2\lfloor \sqrt{s} \rfloor$. 
By a similar argument to the proof of Lemma \ref{affectlem} (this
time we omit the details),
we  have the following result, which is analogous to Lemma 2.3
 of \cite{BR06}.

\begin{lemm} \label{edense}
Let $\rho \geq 1$ be constant. Let $R_s$ and $E_{dense}(R_s)$ be described as above. Let $r = 2\lfloor \sqrt{s} \rfloor$. Then
 $P[E_{dense}(R_s)] \to 1$  as $s \to \infty$.
 Also  $E_{dense}(R_s)$ depends only on the arrival times
at  sites in $R_s[r]$.
% and
% if it occurs then the state of every site in $R_s$ is fixed.
\end{lemm}
%
%
%{\bf Proof.}
% Any path from $R_s$ to the edge of $R_s[r]$ must start on the edge of $R_s$ and have length at least $r+1$. The number of possible paths of length $r+1$ from the edge of $R_s$ is at most $2 (s + \lfloor \rho s \rfloor ) 3^r := Ks 3^r$. The probability that a change can propagate along a path of length $r+1$ is at most $\frac{1}{[r/2]!}$. Any path from the edge of $R_s$ to the edge of $R_s[r]$ must start off with one of these paths of length $r+1$ so the probability that $E_{dense}(R_s)$ does not occur is at most $\frac{Ks 3^r}{[r/2]!}$ which goes to $0$ as $s$ goes to $\infty$. If the event $E_{dense}(R_s)$ occurs then nothing outside $R_s[r]$ can affect the states of the sites inside $R_s$.
%
%\hfill $\Box$

To prove Proposition \ref{RSWlem}, assume
 for a contradiction that it does not hold and fix a value of
 $\lambda$ where it fails. Then
$\liminf_{s \to \infty} f_\lambda^*(1,s)>0$ and for
some $\rho >1$ we have
$\lim_{s \to \infty} f_\lambda^*(\rho,s) = 0$.
But then, as in (4.4)  of \cite{BR06},
 for any $\epsilon > 0$ we have $f^*_{\lambda} (1+\epsilon,s) \rightarrow 0$ as $s$ goes to $\infty$. Throughout the argument let $T_s$ be the strip $[1,s] \times {\mathbb Z}$. The first claim in the proof in \cite{BR06}
can easily be adapted to the integer lattice as follows.

\begin{lemm} \label{infsqRCM}
Let $\epsilon > 0$ be fixed and let
 $\delta := \delta (s) := \frac{\lfloor \epsilon s \rfloor}{s}$.
Let $L$ be the line segment $\{1\} \times [-\delta s, \delta s]$. Then the probability that there is a black path $P$ in $T_s$ starting from $L$ and going outside $S' = [1,s] \times [-(1/2 + 2 \delta) s, (1/2 + 2\delta)s]$ tends to zero as $s \rightarrow \infty$.
\end{lemm}

{\bf Proof.}
By symmetry in the line $[1,s] \times \{0\}$ it suffices to show that the event $E$ that there is a black path $P_1$ lying entirely within $S'$ and connecting some site of $L$ to some site at the top of $S'$ has probability tending to zero.

Let $E_1$ be the event that there is such a path $P_1$ lying entirely in the rectangle $R = [1,s] \times [-s/2, s/2 + 2\delta s]$. If $E$ holds but $E_1$ does not then there is a black crossing the long way of an $s$ by $s +2\delta s + 1$ rectangle which has probability tending to zero. Therefore if suffices to show that $P(E_1) \rightarrow 0$.

Reflecting vertically in the line $y = \delta s$, let $L' := \{1\} \times [\delta s, 3 \delta s]$ be the image of $L$. Let $E_2$ be the event that there is a black path $P_2$ from $L'$ to some point with height $-s/2$. Then by symmetry and by the Harris-FKG inequality the probability that $E_1$ and $E_2$ occur is at least $P(E_1)^2$. But then $P_1$ and $P_2$ must meet and therefore contain a black path crossing $R$ from top to bottom which has probability tending to zero,
so $P(E_1) \rightarrow 0$.
\hfill $\Box$ \\

{\bf Proof of Proposition \ref{RSWlem}.}
The rest of the claims in \cite{BR06}
can be treated similarly by replacing squares in the
 plane with squares in the ${\mathbb Z}^2$ lattice and
 making use of the FKG inequality holding in the RSA model. Then this
 combined with Lemma \ref{edense} and the fact that $r = 2\lfloor \sqrt{s}
 \rfloor$ is $o(s)$ completes the proof of Proposition  \ref{RSWlem}.
\hfill{$\Box$} \\

For $n \in \N$ we define the boxes
\bea
B(2n+1) := [-n,n] \times [-n,n];
~~~
B(2n) := [-n,n-1] \times [-n,n-1].
\lbl{boxdef}
\eea 

%\section{}

{\bf Proof of Theorem \ref{thm1}.}
By Proposition \ref{keyprop}, there exists $\mu<1$ such that
(\ref{eqstar}) holds.  
Defining $ \delta := (1/3) \limsup_{s \rightarrow \infty} f^*_{\mu}(4,s) $,
 we have by (\ref{eqstar})  
and Proposition \ref{RSWlem} that 
 $\delta >0$.  

Thus we can find infinitely many even $n$ such that the probability of a 
black crossing (including diagonals) the long way of a $4n$ by $n$ rectangle
 is at least $2 \delta$. With any such even $n$ we can find an odd $m$ bigger
 than $n$ such that a crossing the long way of a $4n$ by $n$ rectangle
 means that there is a crossing the long way of a $3m$ by $m$ rectangle.
 Therefore we can find infinitely many odd $n$ such that there is a crossing
 of a $3n$ by $n$ rectangle with probability at least $2 \delta $. 
Then for odd $n$, using the Harris-FKG inequality 
 for this model (see \cite{PS}), the probability of there
being a circuit of the annulus $B(3n) \setminus B(n)$ is at least
 $(2\delta)^4 $.

By Lemma \ref{affectlem}, 
for any $n$ we can find an $m$ depending on $n$ such that
$$
P \left[
\cup_{y \in \Z^2 \cap B(n), z \in \Z^2 \setminus B(m) }
(\{ y ~{\rm affects~}z\} 
\cup \{ z ~{\rm affects~}y\} 
) 
\right]
% the probability of the arrival times outside $B(m)$ affecting the colours of the vertices in $B(3n)$ is less than $l/2$. To see this the 
\leq \delta^4.
$$
%probability that a particular vertex in $B(3n)$ is affected is no more than the sum over all self-avoiding paths to points outside $B(m)$ of the probability that the arrival times along this
% path are all in descending order, so that the change can propagate along it. 
%
%
% Therefore the probability that any vertex in $B(3n)$ is affected is less than
 %$\frac{9*n^2*4*3^m}{m!}$ 
% $9(n^2)4(3^m)/\lfloor m/2 \rfloor!$, 
%which is less than $l/4$ for $m$ large enough.
%
% In similar fashion we note that for any $m$ there exists $M$ such that the probability of arrival times inside $B(m)$ affecting the states of any vertices outside $B(M)$ is no more than $l/4$. 
% 
% 
Thus, we can build up a sequence
 $m_1 < n_1 < m_2 < n_2 < ...$ such that (i) for any $i \in \N$,
 the annulus $A_i : = B(3n_i) \setminus B(n_i)$ fits inside the annulus
 $A'_i := B(m_{i+1})\setminus B(m_i)$ and (ii)
The probability that there exists 
any vertex inside $A_i$ that is affected 
from outside $A'_i$ is at most $\delta^4$.

Then let $E_i$ be the event that (i) there is a closed circuit around
the origin consisting of sites in the annulus $A_i$
that are black for the process restricted to $A'_i$ and (ii)
no site of $A_i$ is affected by any site 
outside $A'_i$. 
Then for all $i$, $P[E_i] \geq \delta^4$
 and all the events $E_i$ are independent.
 If any one of these events occurs there cannot be an infinite white component in $\Lambda$ containing the origin, so by the Borel-Cantelli lemma
 the probability of an infinite white
component occurring is $0$. Therefore we have 
$$\lambda_c \geq 1/\mu > 1
$$
which completes the proof, subject to proving Proposition \ref{keyprop}.
% if we can find such a $\mu$ satisfying (\ref{eqstar}).
\hfill{$\Box$} \\

\section{Enhancement}
\lbl{secenhance}
\eqco

\begin{figure}[htbp]
\begin{center}
\scalebox{1.7}{\includegraphics[angle = 0]{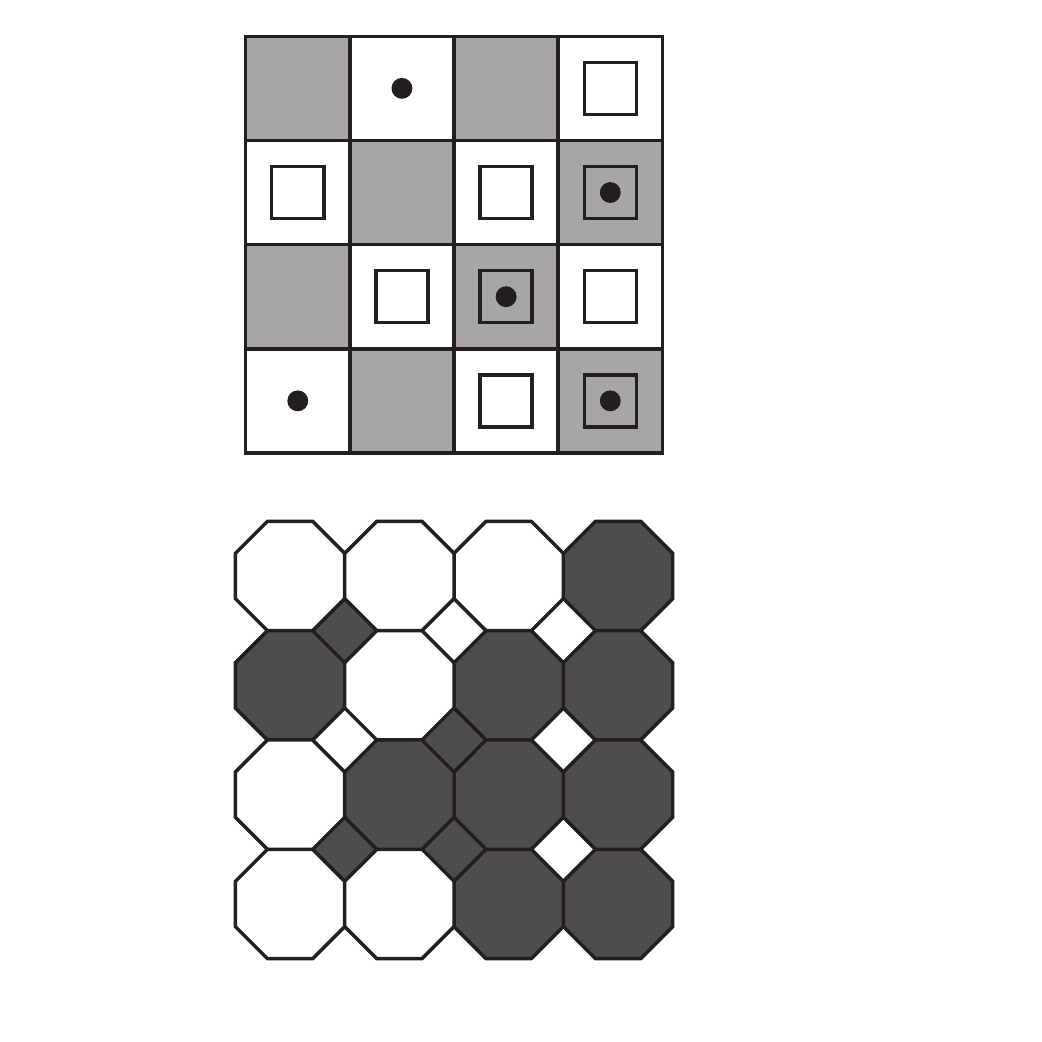}}
\end{center}
\caption{Here is an example of random sequential adsorption and a corresponding
% $K_{8,8,4}$ face
 percolation process on the faces of the $(4,8^2)$ lattice}
\label{fig:RSA.Klat}
\end{figure}

We now define an enhancement that we shall use to interpolate
between the RSA models  on $\Lambda$ and on $\Lambda^*$. Consider
 the infinite $(4,8^2)$ lattice  
(see Figure \ref{fig:RSA.Klat}: we use terminology from \cite{BRbk}, page 155),
with faces  divided into {\em octagons}  and {\em diamonds}.
The octagons are centred at the sites of $\Z^2$, and
the diamonds are centred at the sites $\{z':z \in \Z^2\}$,
where we set $z':= z+ (1/2,1/2)$ (we shall refer to sites
$z', z \in \Z^2$ as {\em diamond sites}).

Now consider a certain  dependent face percolation model on
 the infinite $(4,8^2)$ 
%  $K_{8,8,4}$
 lattice, 
in which each octagon is given the same colour (black or white)
as the corresponding site in the random
 sequential adsorption model, and
 each of the diamonds is black
 with probability $p$ (the {\em enhancement probability})
and white otherwise (independently of
everything else). Thus  $p=0$ is equivalent to
 $\Lambda$ and $p=1$ is equivalent to $\Lambda^*$. 

Placing a vertex at the centre of each face of the $(4,8^2)$
lattice, and taking two vertices to be adjacent
if and only if the corresponding faces 
of the $(4,8^2)$ lattice are adjacent,
we obtain the so-called  {\em centred quadratic lattice}
 (see \cite{Kes}), and we may equivalently view the
dependent face percolation model just described as
 a site percolation model on the centred quadratic lattice.

Let $h(n,\lambda,p)$ denote the probability that there is
 a horizontal black crossing in $\Lambda$ of a
 $2n$ by $2n$ square $B(2n)$ (as defined at (\ref{boxdef}))
 with arrivals rate $\lambda$ on the
 even sites and $1$ on the odd sites and enhancement probability $p$.
In this model we must have either a horizontal crossing or a
 vertical white crossing but not both. Also,
for $(\lambda,p)=(1,0.5)$ the probabilty of both
 these events must be the same by symmetry so the probability
 of a horizontal black crossing is $0.5$.
That is, for any $n$ we have
\bea
h(n,1,0.5)=0.5.
\label{1024a}
\eea

\begin{remk}
\lbl{weakineq}
By (\ref{1024a}) and  monotonicity, we have
$h(n,1,1) \geq 0.5$ and therefore (\ref{eqstar}) holds for $\mu=1$.
Hence, by the argument already given in the proof of Theorem \ref{thm1}
at the end of Section \ref{secdual}, we have $\lambda_c \geq 1$.
The  remainder of this paper is concerned with demonstrating
that this inequality is strict.
\end{remk}

To each diamond $x', x \in \Z^2,$
we assign a uniform random variable $T_{x'}$ 
(the {\em enhancement variable}). 
Then $x'$ is black if $T_{x' }< p$ and white otherwise.
 We then introduce the idea of a site being {\em pivotal}.
 Let $H_n$ be the event that we have a horizontal crossing of $B(2n)$ in the
 enhanced model on $\Lambda$. Then we say that an even site $x$ is 
$1$-pivotal if making the arrival time $t_x$ equal to the first
arrival of the Poisson process at $x$ 
%$0$
 means that $H_n$ occurs but making $t_x $ equal to the
second arrival time of the Poisson process at $x$ 
%= \infty$
 means it does not.
 We say that a diamond $x'$ is $2$-pivotal if making
 $T_{x'} = 0$ means $H_n$ occurs but if $T_{x'} = 1$ then it does not.

For $x \in \Z^2$,
let $P_1(n,\lambda,x)$ be the probabilty that site
$x$ is $1$-pivotal, and
let $P_2(n,\lambda,x)$ be the probabilty that site
$x'$ is $2$-pivotal.
 We have the following proposition (a variant of the Margulis-Russo formula).

\begin{prop} 
\label{Russoprop}
It is the case that
\bea
\frac{\partial h(n,\lambda,p)}{\partial \lambda} = (1/\lambda)
 \sum_{x \in \Z^2: x ~ {\rm even}} P_1 (n,\lambda,x)
\label{1011b}
\eea
and
\bea
\frac{\partial h(n,\lambda,p)}{\partial p} = \sum_{x\in \Z^2} P_2 (n,\lambda,x)
\label{Russo2}
\eea
\end{prop}
{\bf Proof.}
Fix $n$ and $p$.
Enumerate the even sites of $\Z^2$ in some manner as
$x_1,x_2,\ldots$. Given $k \in \N$ and given
$\lambda_1 >0,\lambda_2 >0$, let $E_k(\lambda_1,\lambda_2)$
be the event that $H_n$ occurs when we use a  Poisson arrivals
process of rate $\lambda_1$ at sites $x_1,\ldots,x_{k-1}$
and of rate $\lambda_2$ at sites $x_k,x_{k+1},x_{k+2},\ldots$.
%and let $g_k(\lambda_1,\lambda_2)$
Let $\epsilon >0$. 
For $x \in \Z^2$ 
let $A(x)$ be the event that site $x$ affects some site
in $B_n$.
Since  the probability there is an path with increasing arrival times
starting at $x$ and ending at $B_n$ decays at least exponentially
in the distance from $x$ to $B_n$ (see the proof of
Lemma \ref{affectlem}),
 the sum  $\sum_{k \geq 1} P[A(x_k)]$
converges. Hence
by the first Borel-Cantelli lemma,
\bean
0 \leq 
P[E_k(\lambda,\lambda+\epsilon) ]
- h(n,\lambda,p)  
\leq 
P \left[ \cup_{j=k}^\infty A(x_j) \right]
\\
\to 0 ~~~{\rm as} ~ k \to \infty.
\eean
Hence,
\bea
h(n,\lambda+\eps,p) -
h(n,\lambda,p) = 
P[E_1(\lambda,\lambda + \eps)]   
- \lim_{k \to \infty} 
 P[E_k(\lambda,\lambda + \eps)]   
\nonumber \\
= \sum_{k=1}^\infty
 P[E_{k}(\lambda,\lambda + \eps) \setminus   
E_{k+1}(\lambda,\lambda + \eps)].   
\lbl{1021b}
\eea
Here we are assuming the Poisson processes of rate
$\lambda$ and $\lambda + \eps$ at
$x_k$ are coupled in the usual way, i.e.
with the the  $(\lambda + \eps)$-process
decomposed into two independent processes
of rate $\lambda$ and $\eps $ respectively.

Event
 $E_{k}(\lambda,\lambda + \eps) \setminus   
E_{k+1}(\lambda,\lambda + \eps)$ occurs
if and only if (i) the first arrival time $T_1$ of
the $(\lambda + \eps)$-process at $x_k$   
comes from the $\eps$-process, and (ii) the
crossing of $B(2n)$ occurs if we use
the arrival  time $T_1$ at $x_k$, but not if
we use the arrival time $T_1 + T_2$,
where $T_2$ is the time
from $T_1$ to the next arrival of the $\lambda$-process at $x_k$.
Note that $T_2$  is exponential with parameter $\lambda$,
independent of $T_1$ and the type of the arrival at time $T_1$.
Therefore, 
\bea
 P[E_{k}(\lambda,\lambda + \eps) \setminus   
E_{k+1}(\lambda,\lambda + \eps)]
= (\eps/(\lambda + \eps)) P[F_k(\lambda,\lambda+\eps)]
\lbl{1021a}
\eea
where $F_k$ denotes the event
that the crossing of $B_n$ occurs
if we use the first arrival at $x_k$
but not if we use the second arrival
at $x_k$, and our arrivals  processes
are Poisson rate $\lambda$ at sites $x_j,j <k$,
 and Poisson rate $\lambda + \eps$ at sites
$x_i, i >k$, and our first arrival
at $x_k$ is exponential rate $\lambda + \eps$
but the time from the first arrival to the  second arrival at $x_k$
is exponential rate $\lambda$.

Coupling events
$F_k(\lambda, \lambda+ \eps) $
and
$F_k(\lambda, \lambda)$, we have for any integer $K >n$
  that
$P[F_k(\lambda, \lambda+ \eps) \setminus F_k(\lambda, \lambda)]$
is bounded by the sum of the probability that
there is some 
site inside $B(2n)$
that is affected from outside $B(2K)$, 
 and the probability that there exists
 some site $x_j$ inside $B(2K)$  such that
the first  arrival for the $(\lambda + \eps)$-process
at $x_j$
comes from the $\eps$-process at that site.  
For any fixed $K$
the second of these probabilities tends to zero as $\eps \downarrow 0$, 
while the first probability is small for large $K$, uniformly in $\eps$.
Hence by (\ref{1021a}),
\bean
\lim_{\eps \downarrow 0}  \eps^{-1}
 P[E_{k}(\lambda,\lambda + \eps) \setminus   
E_{k+1}(\lambda,\lambda + \eps)]
= \lambda^{-1}
P[F_k(\lambda, \lambda) ] 
%nonumber \\
= \lambda^{-1}
 P_1(n,\lambda,x_k).
\eean
Moreover,  $P[F_k(\lambda,\lambda+\eps)]$ is
bounded by the probability that  there are
 increasing arrival times
along
some path from $x_k$ to $B_n$, which is
bounded by a summable function of $k$ uniformly
in $\eps$. Therefore by (\ref{1021b}), (\ref{1021a})
 and dominated convergence we have
\bea
\frac{\partial ^+ h}{\partial \lambda} = \lim_{\eps \downarrow 0}
 \frac{h(n,\lambda + \eps ,p) - h(n,\lambda ,p)}{\eps}
= \lambda^{-1} \sum_{k=1}^\infty  P_1(n,\lambda,x_k).
\eea
By a  similar argument (we omit details), one can obtain the same expression
for the left derivative
$\frac{\partial ^- h}{\partial \lambda} $. Therefore
(\ref{1011b}) is proven.

The proof for the second part (\ref{Russo2})
 is similar.
\hfill{$\Box$}

\section{Comparison of pivotal probabilities}
\label{secpivcomp}
\allco
The following proposition is a key step in the proof of Theorem \ref{thm1}.
\begin{prop}
There exists a constant $K_1 \in (0,\infty)$
such that for all $n$,
 all $(\lambda,p) \in [0.5,1.5] \times [0.2,0.8]$,
and any even site 
$y$ in $B(2n)$,
 there exists an adjacent diamond site $\ty'$ such that
$$
 P_1(n,\lambda,y)
\leq 
K_1 P_2(n,\lambda,\ty). 
$$
\label{propP1P2}
\end{prop}
The rest of this section is devoted to proving 
Proposition \ref{propP1P2}.
The argument  is quite lengthy and we
 divide it into stages.

Fix $y \in \Z^2$. For $r,s \in \N$ with $s >r$, let
$C_r$ be the square of side $2r+1$ centred at $y$,
and let $A_{r,s} := C_s \setminus C_r$.

We shall consider  a coupling of RSA processes. Let
 $S_x$ be the arrival times and
enhancement variables in one process (so if $x \in \Z^2$ then
$S_x$ is exponentially
distributed but $S_{x'}$ is a uniformly
distributed enhancement variable).
Let $T_x$ be the arrival times and enhancement variables
 in another independent process.  
 Given $r,s \in \N$ with $s \geq r$,
we use these to
create a third process of arrival times and enhancement variables
 $U_x^{(r,s)}$, as follows.
Put
\bea
 %Let $U_x = S_x$ and $U'_x = S'_x$ for all $x$ outside
U_x^{(r,s)} :=  \left\{ \begin{array}{ll} S_x, &  x \notin C_{s} \\
  {\cal B}_x S_x + (1- {\cal B}_x) T_x, & x \in A_{r,s} \\
T_x, & x \in C_{r}
\end{array}
\right.
\lbl{1014d}
\eea
 where the $ {\cal B}_x$ are independent Bernoulli variables
 with parameter $0.5$. 

The next lemma establishes a sort of contitional independence
between the occupancy status, in the $U_x^{(r,s)}$
process, of sites inside $C_r$
and of sites outside $C_s$, conditional on the occurrence
of a certain event associated with sites in the
annulus $A_{r-2,s}$.

For $x \in \Z^2$,
 define $I_S(x)$ to be $1$ if site $x$ is occupied and $0$ if 
it is blocked in the $S_x$ process. 
Define the following sets of sites:
\bea
%M^{(r,s)}: = \{ x \in A_{r+1,s+1} \cap\Z^2 : I_S (x)=1 \};
M^{(r,s)}: = \{ x \in A_{r,s} \cap\Z^2 : I_S (x)=1 \};
~~~~~~
N^{(r,s)} := A_{r,s} \setminus M^{(r,s)};
\lbl{1024b}
\eea
%\[
%N^{(r,s)} := \{ x \in A_{r,s} : {\rm  dist} (x,M ) = 1 \};
%\]
\[
M^{(r,s)}_1: = \{ x \in M^{(r,s)} : S_x \leq 1 \}; ~~~ 
M^{(r,s)}_2 := M^{(r,s)} \setminus M^{(r,s)}_1 ;
\]
\[
N^{(r,s)}_1 := \{ x \in N^{(r,s)} : S_x \leq 1 \};~~~
N^{(r,s)}_2 := N^{(r,s)} \setminus N^{(r,s)}_1.
\]
Define the event
\bea
%E_2 := 
E_1^{(r,s)} :=
  \cap_{x \in M^{(r,s)}_1 \cup N^{(r,s)}_2}
%\{  U_x = S_x \}
\{  {\cal B}_x = 1 \}
\cap 
  \cap_{x \in M^{(r,s)}_2 \cup N^{(r,s)}_1}
%\cap_{x \in A_{2r+6,2r+9} \setminus (M_1^{(r,s)} \cup N^{(r,s)}_2) }
\{ {\cal B}_x = 0 \} 
\nonumber
\\
%; \eeaa
%\beaa
%E_3 :=  E_3(y,r) :=
\cap  \cap_{x \in M^{(r,s)}_2} \{ T_x \leq 1  \} \cap 
\cap_{x \in N^{(r,s)}_1} \{ T_x >1 \}
\nonumber
\\
\cap \cap_{x \in M^{(r-2,r)}} \{ T_x  \leq 1 \}
\cap \cap_{x \in N^{(r-2,r)}} \{ T_x  > 1 \}
.
\lbl{E1def}
\eea
\begin{lemm}
\lbl{Rlem}
Suppose $r,s \in \N$ with $r \geq 3$ and $s \geq r+3$.
If $E_1^{(r,s)}$ occurs then the state of
 all sites in $\Z^2 \setminus C_{r}$ is the same in
 the $U_x$ process as in the
 $S_x$ process.
\end{lemm}
{\bf Proof.}
Assume event 
 $E_1^{(r,s)}$ occurs.
We start off with all the arrival times in the $S_x$ process. 
Then we change 
the arrival times in $M^{(r,s)}_2$ one by one.
 Each time we are making the arrival
 time at an occupied site earlier, so we cannot change the state of any sites.
 Then we change the arrival times in $N^{(r,s)}_1$ one by one. Each time we are
 making the arrival time at a blocked site later so we cannot change the state
 of any site. We then have our $U_x$ process on $\Z^2 \setminus C_r$.

Now we change the arrival times for  the sites inside $C_r$. 
Every site  $x \in M^{(r-1,s-1)}$ has $U^{(r,s)}_x \leq 1$ and
has all its neighbours $z$ with $U^{(r,s)}_z >1$, so is occupied
in the $U^{(r,s)}$-process.
Also, every site $z \in N^{(r,s-2)}$ has
  $U^{(r,s)}_z > 1$ and
has at least one occupied neighbour $x $ with
  $U^{(r,s)}_x \leq 1$, so is vacant.

Thus when we change the arrival times for the sites inside $C_r$, the
states of sites in $A_{r,s-2}$ do not change
and therefore the states of sites in $\Z^2 \setminus C_{s-2}$ also
do not change.

  Hence, whatever arrival times we have on
$C_{r-2}$, the states of the sites of $\Z^2 \setminus C_r$ do not change,
 so they are the same in the $U_x^{(r,s)}$ process as
 in the $S_x$ process. 
\hfill{$\Box$} \\

We aim to prove Proposition \ref{propP1P2}, so let us
assume $y \in B(2n)$ and $y$ is an even site.
 Let $\ty'$ be the first diamond site adjacent to
 $y$ that is contained in $B(2n)$
 working clockwise from the top right (so $\ty=y$ if $y$ is 
in the interior of $B(2n)$). Let 
$D_r$ be the diamond of sites that are at $\ell_1$ distance $r$ or less
 from $y$.

We shall say that
$y$ is $(1,r)$-$pivotal$ for event $H_n$
if changing  $t_y $ from the second Poisson arrival
time to the first arrival time, and changing any affected sites within $r$ steps of $y$, means that $H_n$ occurs but
 changing only sites within $r-1$ steps of $y$ means $H_n$ does not 
occur (by this we mean changing the $4$ sites adjacent to $y$ as appropriate as the first step 
then changing any sites adjacent to these as appropriate as the second step and so on).  Define $P_{1,r,n}(y)$ to
 be the probability that  $y$ is $(1,r)$-pivotal for $H_n$.

Given $n$ and $y$, define event $E(r)$,
for $r \in \N$, as follows. 
First suppose that
 $r \leq n/5$ and
the left and right endpoints of
 $D_{r+7} $ lie in
$ B(2n)$.
Then let $E(r)$ be the event
 that we have black paths in $B(2n)$ from each side of $B(2n)$ up to 
$\Z^2 \cap C_{r+7} $ 
but no black path from one side of $B(2n)$ to the other avoiding
$C_{r+7}$. Here we are using the second arrival time at $y$. 

If $r \leq n/5$, and the left (respectively, right)
  endpoint of $D_{r+7} $ lies outside $ B(2n)$, then
let $E(r)$ be the event
 that we have a black path in $B(2n)$ from the right (respectively, left)
  side of $B(2n)$ up to 
$\Z^2 \cap C_{r+7} $, but no black path in $B(2n)$
from one side of $B(2n)$ to the other avoiding $C_{r+7}$.

If $r >  n/5$ 
then we define $E(r)$ to be the whole sample
space, so that $P[E(r)]=1$.

\begin{lemm}
\label{newElem}
There exists a constant $K_2 \in (0,\infty) $
 such that for all $n,r \in \N$ and all
even $y \in \Z^2 $, we have
\bea
P_{1,r,n}(y) \leq \frac{  K_2^r  P [E(r) ] }{ \lfloor r/2 \rfloor! },
~~~~~ r \geq 20; 
\lbl{1014b}
\\
P_{1,r,n}(y) \leq    P [E(r) ], 
~~~~~ r \leq 20.
\lbl{1111a}
\eea
\end{lemm}

 Suppose  $y$ is $(1,r)$-pivotal.
  Then, after changing all sites affected up
to  $r$ steps from $y$ when we set $t_y$ to be the first
arrival time rather than the second arrival time,
 we obtain a black crossing of
 $B(2n)$. Any such crossing  path must include at least
one site in 
 $D_r$ (otherwise $y$ would not be $r$-pivotal). 
Therefore event $E(r)$  occurs.
Since $P[E(r)]$ is nondecreasing in $r$, 
 this immediately
gives us (\ref{1111a}).

Now suppose $r \geq 20$.
Let $F(r)$ be the event that there is a self-avoiding path in
$\Z^2$ from $y$ of
 length $r$, namely $y_1, y_2, y_3, \ldots, y_r$, such that 
$t_{y_1} < t_{y_2} < \cdots < t_{y_r}$. 
If  $y$ is $(1,r)$-pivotal
 then $F(r)$ must occur, and hence
% (and so must  $E$  by Lemma \ref{lemE}).
%Hence
\bea
\lbl{1024c}
P_{1,r,n}(y) \leq P[E(r) \cap F (r)].
\eea
Also, as in the proof of Lemma \ref{affectlem} we have
\bea
P[F(r)] \leq
 \frac{ 4( 3^{r-1} )}{ \lfloor r/2 \rfloor ! } 
\lbl{1014a}
\eea
and $F(r)$ depends only on the arrival times inside $D_{r}$.
However, it is not independent of $E(r)$.

We now consider the independent families of arrival 
times $(S_x)$ and $( T_x)$, and a coupled
arrival time process $ U^{(r+2, r+5)}_x$
as defined by (\ref{1014d}).

Let $E^S$, respectively $E^U$, be the event that $E(r)$ occurs based
on the $S_x$ process, respectively the   
 $U^{(r+2, r+5)}_x$ process.
Let $F^S$, respectively $F^U$ be the event that $F(r)$ occurs based
on the $S_x$ process, respectively the   
 $U^{(r+2, r+5)}_x$ process.
Then, defining event 
$A:= E_1^{(r+2,r+5)}$ as given by
 (\ref{E1def}), we have
from Lemma \ref{Rlem} the event identity  $E^S  \cap A = E^U \cap A$. 
Hence,
\bean
P[E^S  \cap F^S] P[A| E^S \cap F^S]
= P[E^S \cap F^S \cap A ] 
\\
= P[E^U \cap F^S \cap A ] 
\\
\leq P[E^U \cap F^S  ] = P[E^U] P[F^S].
\eean 
Also, there is a constant $K_3$ such that
$$
P[A|E^S \cap F^S]
%=P[A|E^S]
 \geq K_3^{-r}.
$$
Combining these 
inequalities
and using the fact that $P[E^U ]= P[E^S]$ yields 
$$
P[E^S  \cap F^S]  \leq K_3^r P[E^S] P[ F^S]
$$
and combined with (\ref{1024c}) and (\ref{1014a}) this gives us the
desired result (\ref{1014b}).
\hfill {$\Box$} \\

Now, given $r \geq 20$,
 we consider for a while the process $U_x:= U_x^{(2r+6,2r+10)}$
as defined by (\ref{1014d}).
%
% and the evek region  
%$R:= R^{(2r+6,2r+9)}$
%as defined by (\ref{Rdef}).
%Let $C_s$ be the box of side $2s+1$ centred on $y$, and
%for $r_1 < r_2$ let $A_{r_1,r_2}$ be the region $C_{r_2} \setminus C_{r_1}$.
%
 Let $G_r$ be the octagonal region 
$C_{2r} \cap D_{4r-10}$, 
 a sort of truncated square.
Note that each of the inner diagonal boundaries of $G_r$ 
consists of  odd sites and is of length 10.
The exact length is not important; we just 
need a reasonably large separation between
each corner of the octagon $G_r$. 
Let $G_r^-$ be the slightly smaller octagonal region 
$C_{2r-4} \cap D_{4r-14}$.

\begin{lemm}
\lbl{Qprop}
There exists a constant $\beta \in (0,\infty)$
with the following property.
Given $r \geq 20$,
if the event $E(r)$ occurs in the $S_x$ process,
then there exists
a stable  set $Q_1 \subset G_r \cap \Z^2$ having
no element  adjacent to the occupied $\Z^2$ sites of the $S_x$
process outside
$G_r$, and
disjoint sets $Q_2, Q_3$ of diamond sites 
inside $G_{r}$,
such that (i) each of $Q_1, Q_2, Q_3$ has at most $\beta r$
elements, and (ii) if, 
in the $U_x$ process, all the sites in $Q_1$ are occupied, 
   all diamonds in $Q_2$ are black,  all the diamonds in $Q_3$ are white,
and all sites in $C_{2r+6} \setminus Q_r$ are in the
same state as for the $S_x$ process,
 then $\ty'$ is $2$-pivotal
for the $U_x$  process.
\end{lemm}
{\bf Proof.}
First suppose $r \leq  n/5$.
 Since $E(r)$ occurs, there must be disjoint black paths
 in the $S_x$ process
up to $\Z^2 \cap C_{r+7} $ from each
 side of $B(2n)$.
 The strategy of the proof is to extend these paths in towards 
$y$ while keeping them disjoint in order to make $\tilde{y}'$ $2$-pivotal.

 For now we assume $C_{2r}$ (and hence $G_r$) is contained in $B(2n)$
(so that $\tilde{y}=y$). 
Let $V$ be the set of black vertices 
(for the $S_x$ proces) in $B(2n) \setminus G_r$
 that are connected to the left hand side of $B(2n)$ by
a black path of the $S_x$ process,
 without using any sites in $G_r$. 
Let $v$ be the first even site
inside $G_r$
 (according to the lexicographic  ordering)
that is occupied (for the $S_x$ process) 
 and connects to $V$ either directly or via blocked odd sites
 adjacent
 to itself and $V$ (and possibly also a black diamond site). 
Let $W$ be the set of black sites 
(for the $S_x$ proces) in $B(2n) \setminus G_r$
that are connected
 to the right hand side of $B(2n)$
by a black path of the $S_x$ process that avoids $G_r$.
Let $w$ be  the first even site that is occupied inside
 $G_r$ and connects to $W$. We now try and build paths from
 $v$ and $w$ in towards $y'$ to make it $2$-pivotal.
 We consider various cases of where $v$ and $w$ are: \\

\begin{figure}%[htbp]
\begin{center}
\scalebox{1.2}{\includegraphics{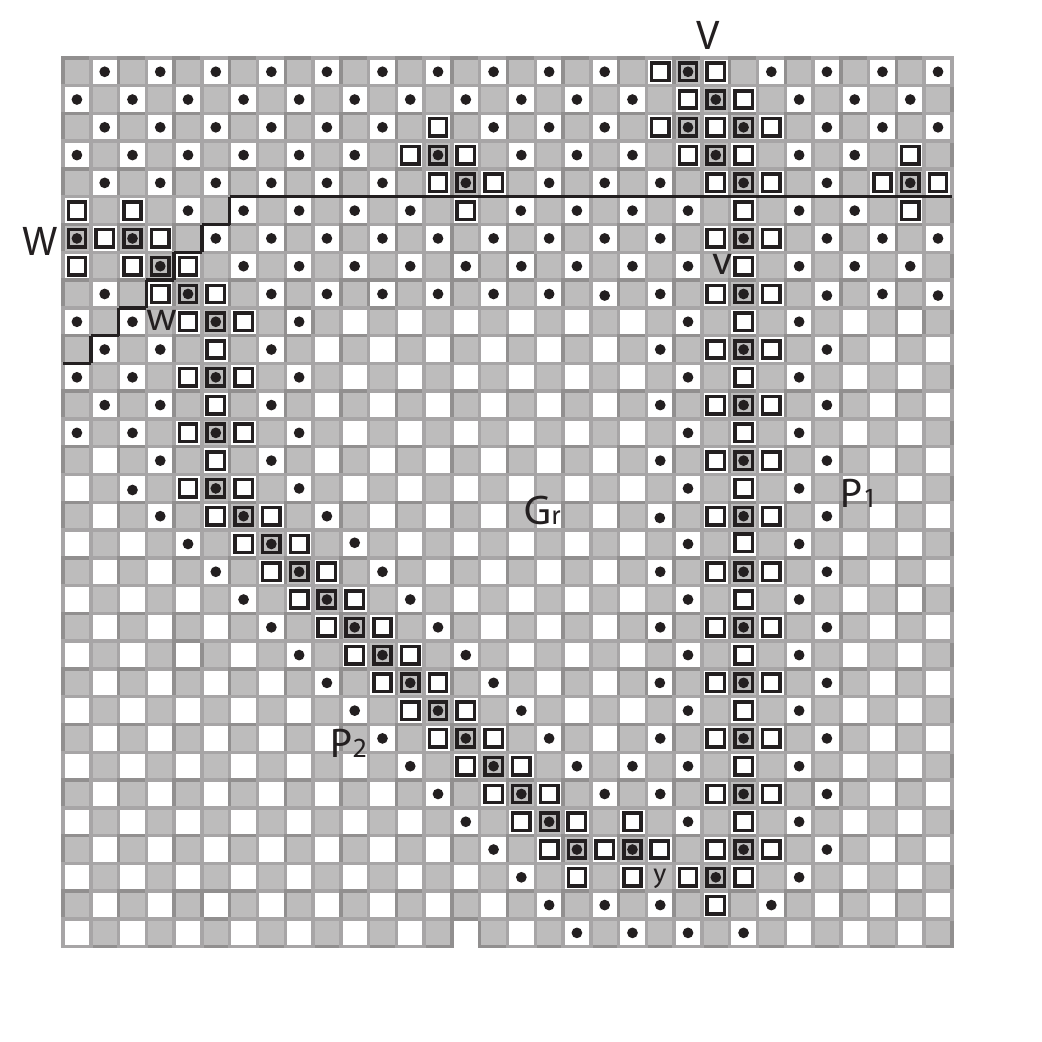}}
\end{center}
\caption{Construction of paths $P_1$, $P_2$ making
$y'$ 2-pivotal.}
\label{paths}
\end{figure}

{\bf Case 1:}
Suppose $v$ and $w$
are well away from each other. In this case we can always make $y'$ 
$2$-pivotal. For example, if $v$ and $w$ are as in Figure \ref{paths},
 we can form disjoint paths $P_1$, $P_2$ of even sites in towards $y$. 
In this and subsequent diagrams, the chequerboard
squares are centred at sites of
$\Z^2$ and are shaded for even sites.
Let $I$ be the set of even sites 
 $\{v , w \} \cup P_1 \cup P_2$. 
Let $J$ be the set of odd sites in 
$G_r \setminus G^-_r$ that are not adjacent to
 any site in $I  $ or to any of the occupied
sites in $C_{2r+6} \setminus G_r$.
Let $J'$ be the set of odd sites in $G_r^-$ that are 
three steps (in $\Z^2$) away from $I$.
 Set $Q_1 := I \cup J \cup J'$.
 %Then 
If the sites in $Q_1$ are occupied for the $T_x$ process,
%arrival  
%times satisfy $T_x < 0.5$ for all the even sites
 %in $I$ and all the odd sites in $J$ and 
 %$T_x > 0.5$ for all other sites in
% $C_{r+10}\setminus R$,
 then  $y$ is $2$-pivotal.
The number of sites in $Q_1$ is bounded by a constant times $r$.
% This all happens with probability at least $(0.2p) ^ {(r+9)^2}$.

\begin{figure}%[htbp]
\begin{center}
\scalebox{1.2}{\includegraphics{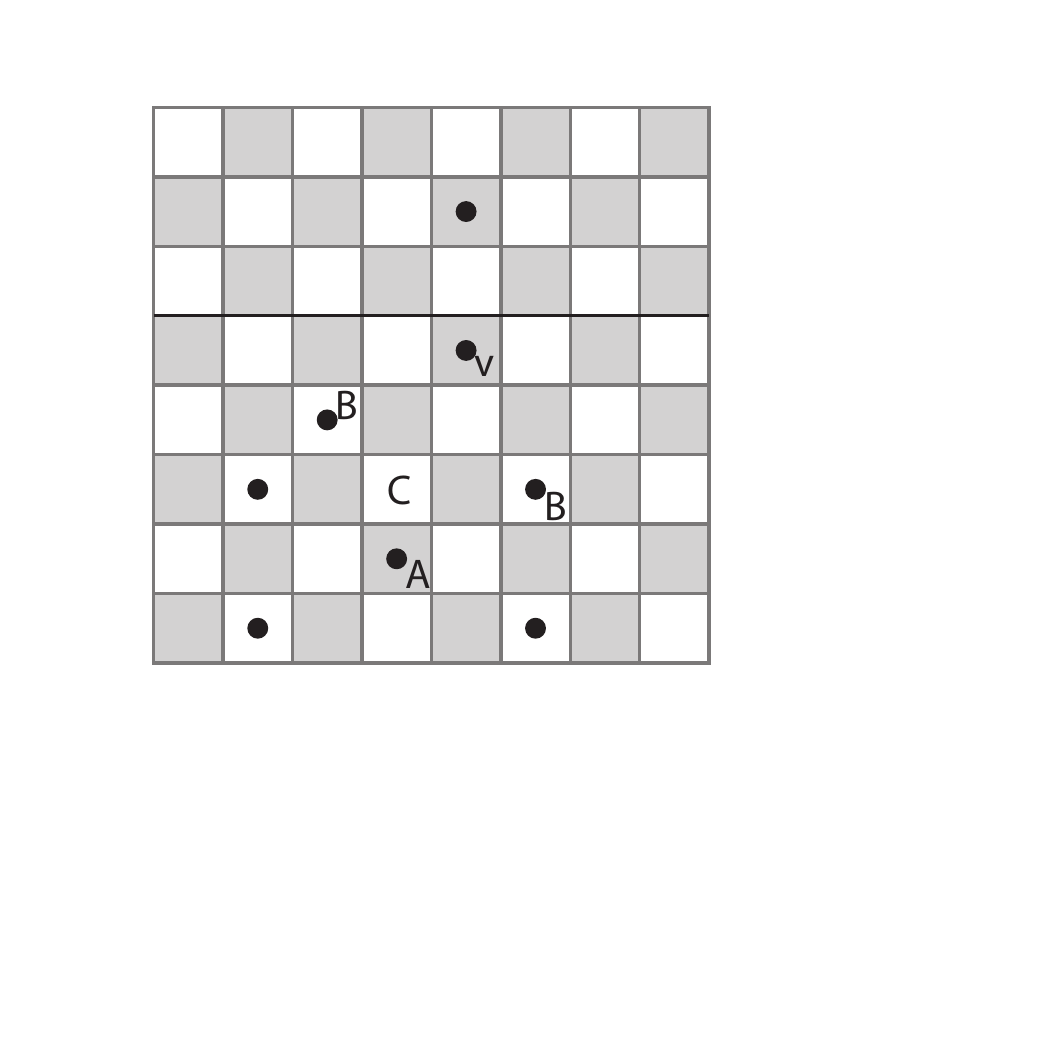}}
\end{center}
\caption{Starting path $P_1$ when  $v$ is on a horizontal edge on the inner
perimeter of $G_r$.}
\label{straight}
\end{figure}

In general, if we have $v$ on a horizontal or vertical edge of
$G_r$, then (see Figure \ref{straight}) 
we can make the even site at position $A$ in relation to $v$ occupied to
 start $P_1$, switch the enhancement on at $C'$ and due to the odd sites at $B$ being occupied
this cannot complete a crossing of $B(2n)$.

\begin{figure}%[htbp]
\begin{center}
\scalebox{1.2}{\includegraphics{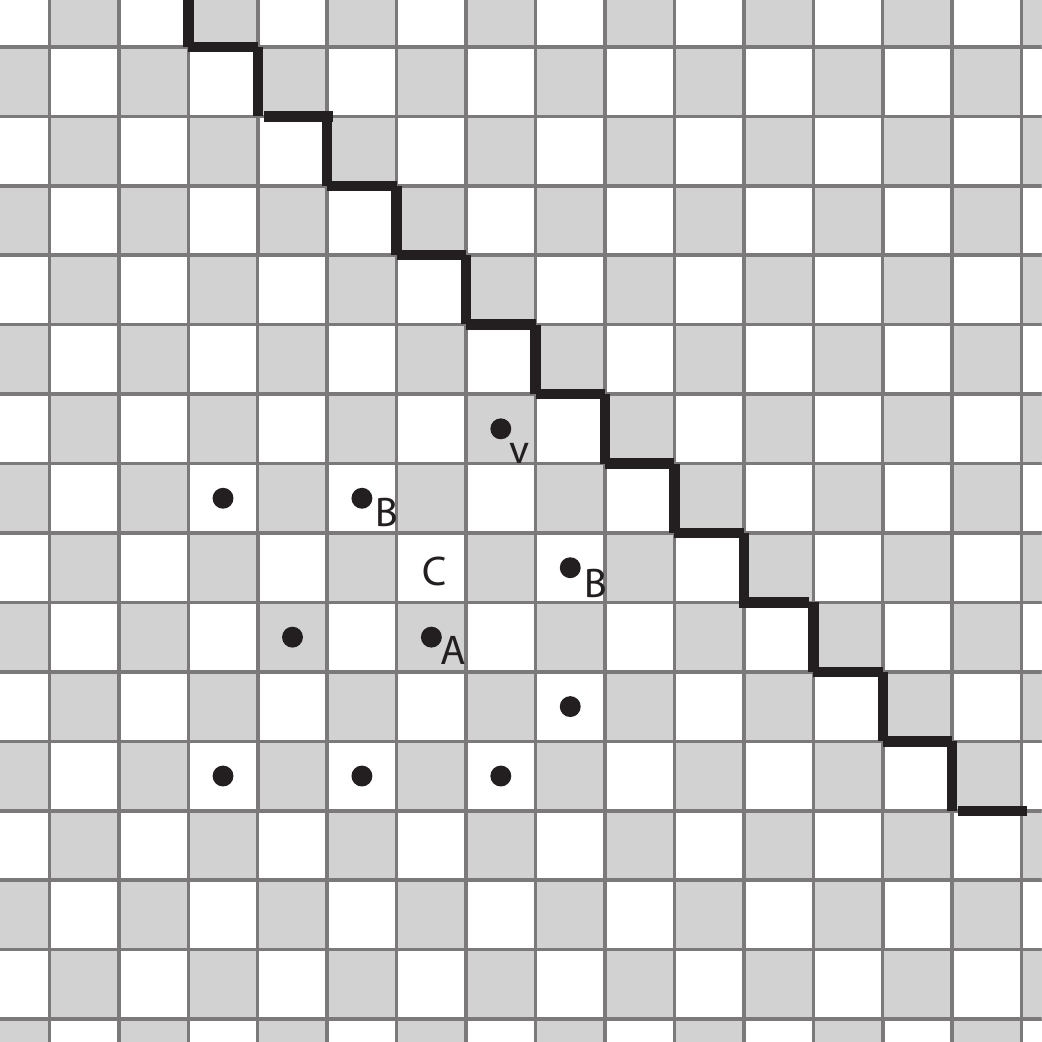}}
\end{center}
\caption{Starting the path $P_1$ when  $v$ lies near a diagonal edge.}
\label{diagonal}
\end{figure}

If $v$ lies beside  a diagonal edge of $G_r$,
 then (see Figure \ref{diagonal})
 we can make the even site at position $A$ in relation to $v$ occupied to 
start $P_1$, switch the enhancement on at $C'$ and due to the odd sites at $B$ being occupied
this cannot complete a crossing of $B(2n)$. \\

\begin{figure}%[htbp]
\begin{center}
\scalebox{1.2}{\includegraphics{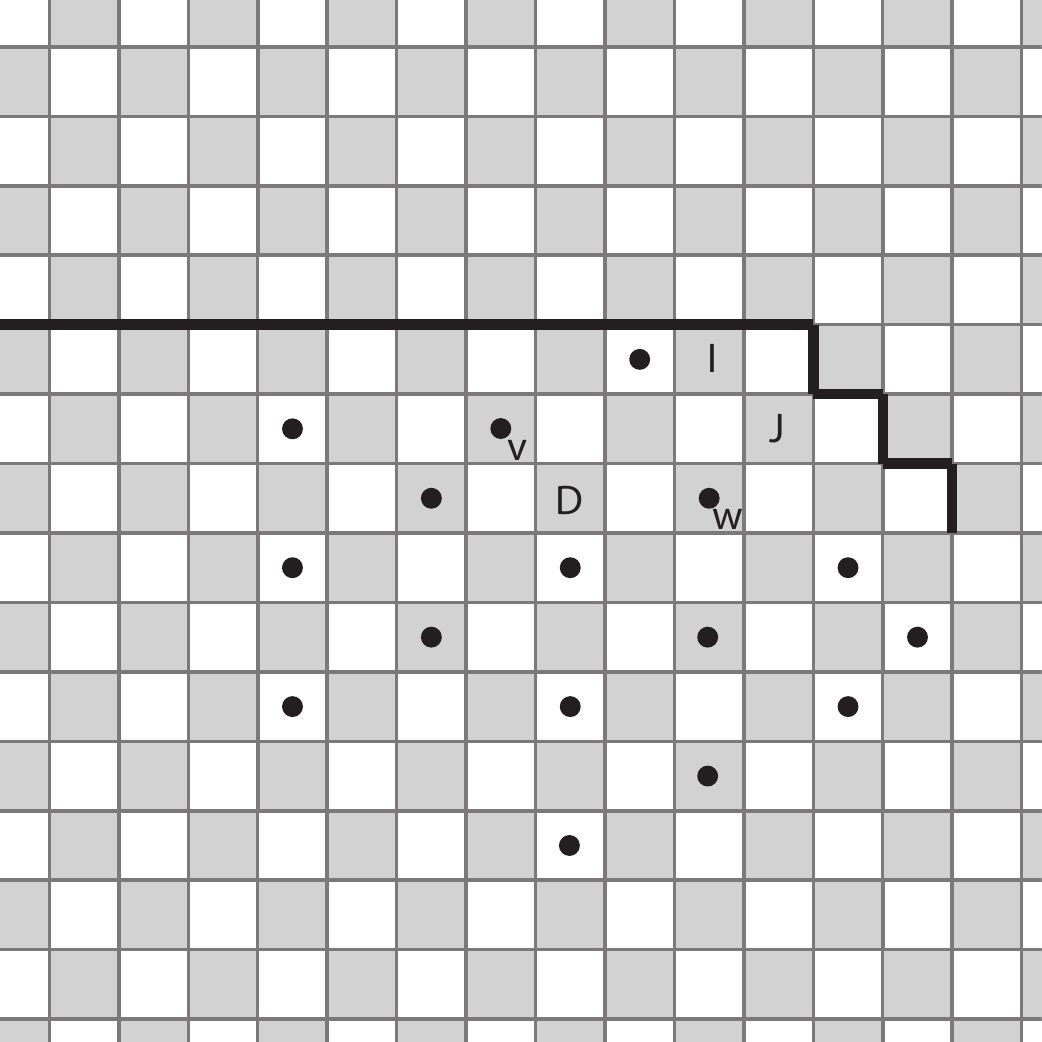}}
\end{center}
\caption{Case 2.}
\label{straight close}
\end{figure}

{\bf Case 2:}
Suppose $v$ and $w$
are near each other but on a straight edge. If their columns are at distance $4$ or more from each other and neither is in position $I$ then there is no 
problem. Their columns cannot be at distance $2$ from each other as then 
$v$ and $w$ would be connected to each other via black sites.
 If they are at distance $3$ then there is no problem
as long as neither $v$ nor $w$ is at position $I$
 (see Figure \ref{straight close}.)
We have the enhancement switched off at $D$ and then extend the paths in towards $y$.  \\

\begin{figure}%[htbp]
\begin{center}
\scalebox{1.2}{\includegraphics{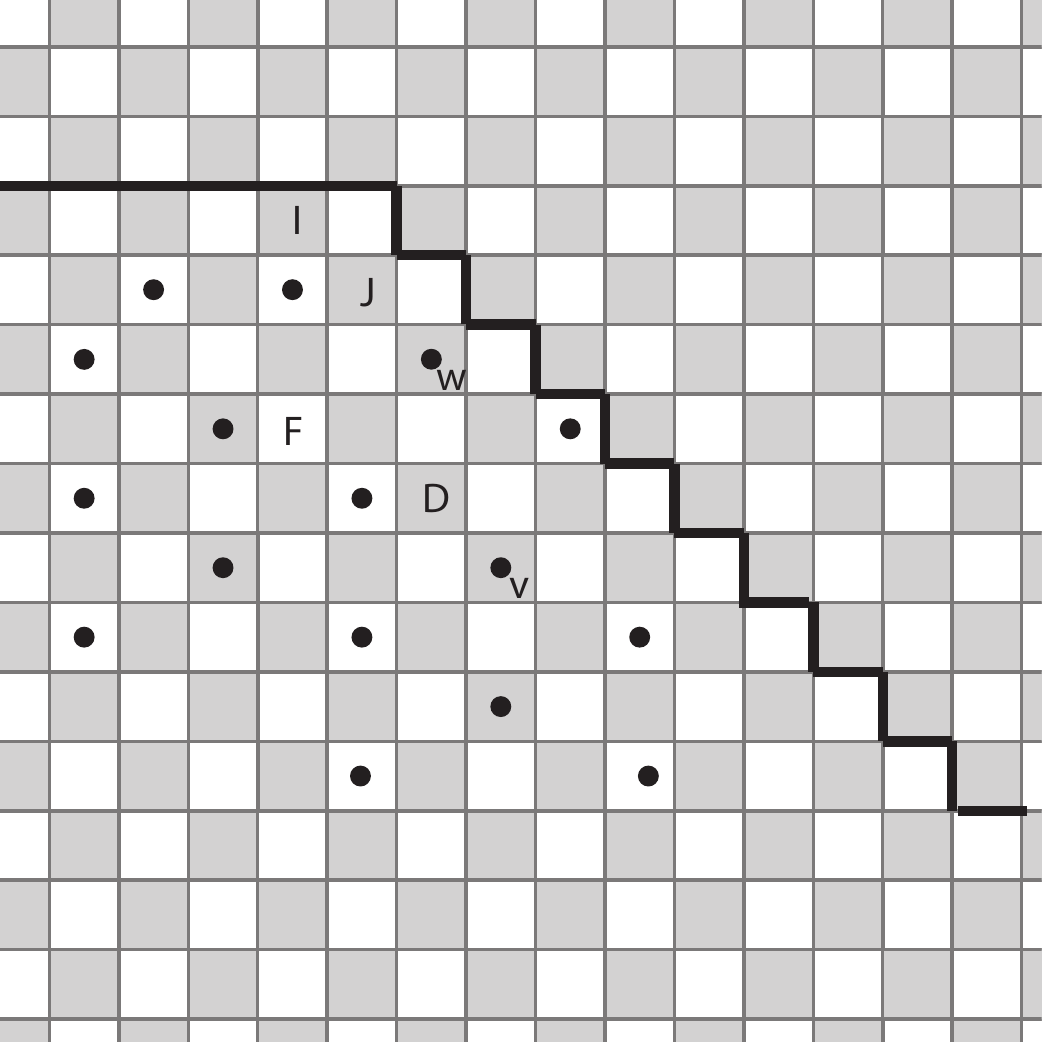}}
\end{center}
\caption{Case 3.}
\label{diagonal close}
\end{figure}

{\bf Case 3:}
Now suppose $v$ and $w$ are near each other on a diagonal edge. If their diagonals are at distance $3$ there is no problem.
 They cannot be at distance $1$ as then they would not be disjoint. If 
they are at distance $2$ and neither is at $J$
(see Figure \ref {diagonal close}),
 there is no problem.
  We have the enhancement switched off at $D$ and switched on at $F$. 
 \\

\begin{figure}%[htbp]
\begin{center}
\scalebox{1.2}{\includegraphics{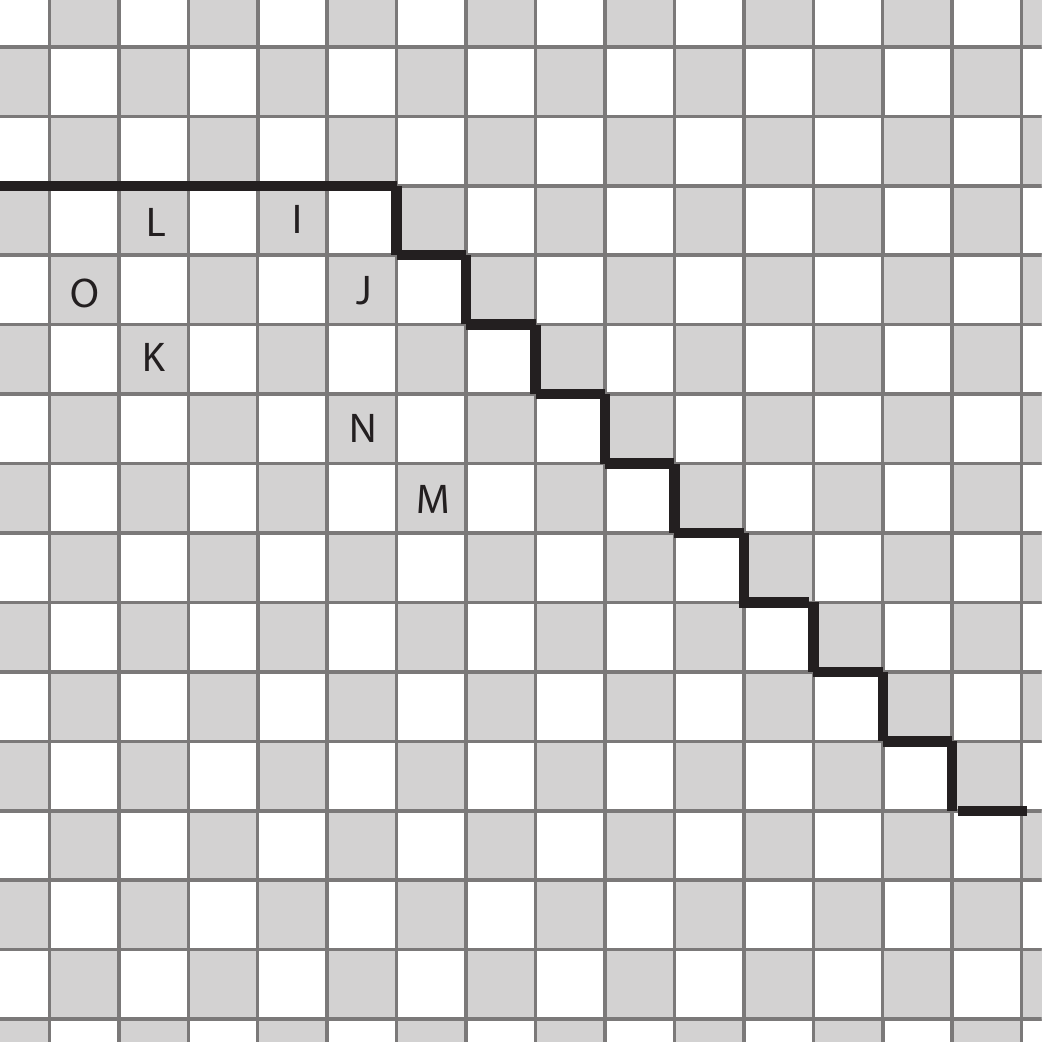}}
\end{center}
\caption{Identifying locations near a corner.}
\label{spread1}
\end{figure}

{\bf Case 4:} Suppose $v$ and $w$ lie
 near to each other but on a corner. We need to consider possible cases when
 $v$ is at $I$ or $J$ (see Figure \ref{spread1}).

(a) $v$ is at $J$. If $w$ is $3$ or more diagonals away then there is no problem. If $w$ is $4$ or more columns away then there is no problem. This just
 leaves three possibilities.

\begin{figure}%[htbp]
\begin{center}
\scalebox{1.2}{\includegraphics{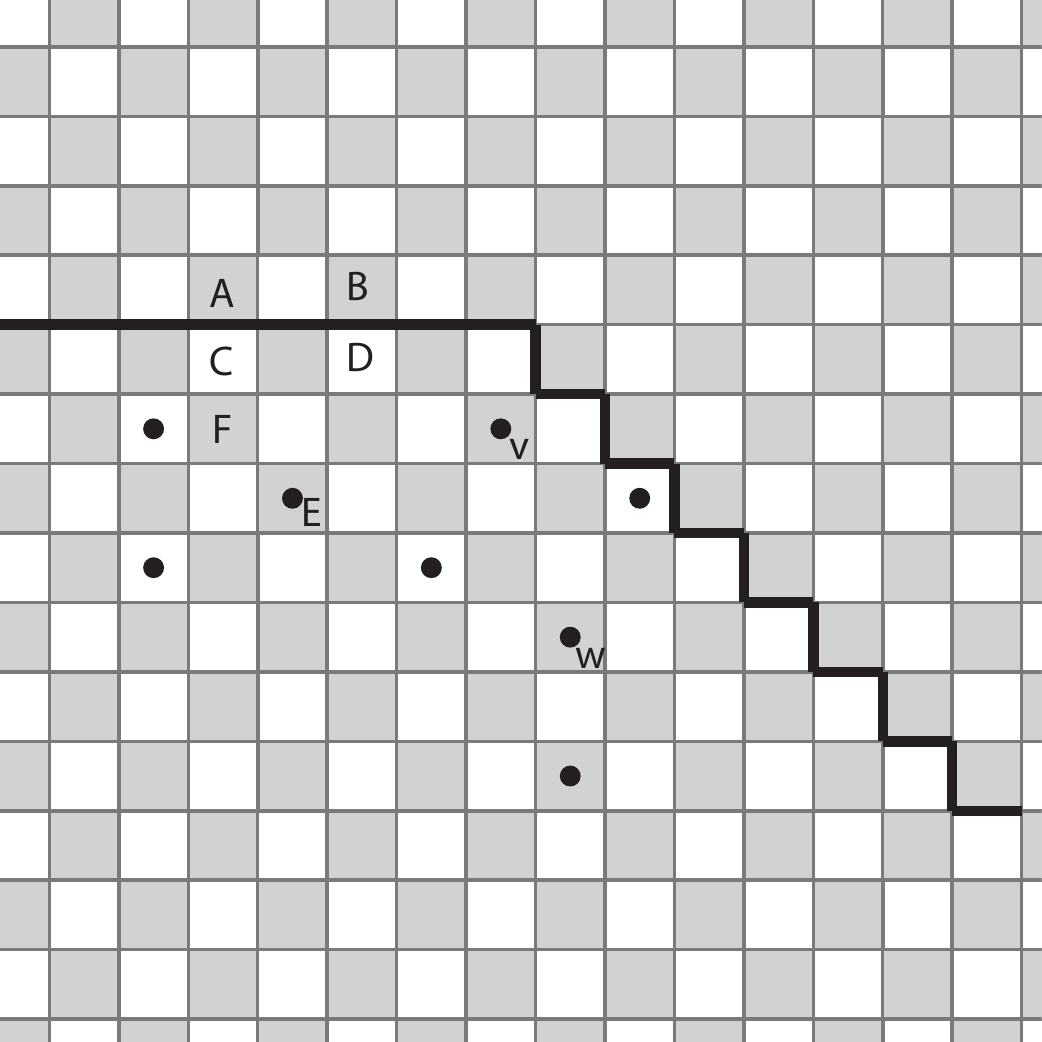}}
\end{center}
\caption{Case 4 (a) (i).}
\label{spread2}
\end{figure}

(i) $w$ is at $M$ (of Figure \ref{spread1}).
Then refer to Figure \ref{spread2}. 
 We can have an occupied even site at $E$, connected
to $v$ via a diamond site.
There is no problem unless there is an occupied even site
 at $A$ that is in $W$. Then we need to have an occupied odd site at $D$ and have the enhancement at  $F'$ switched off.
 We can make $D$ occupied because we know $B$ is unoccupied since
 otherwise it would connect to both $v$ and $W$.

\begin{figure}%[htbp]
\begin{center}
\scalebox{1.2}{\includegraphics{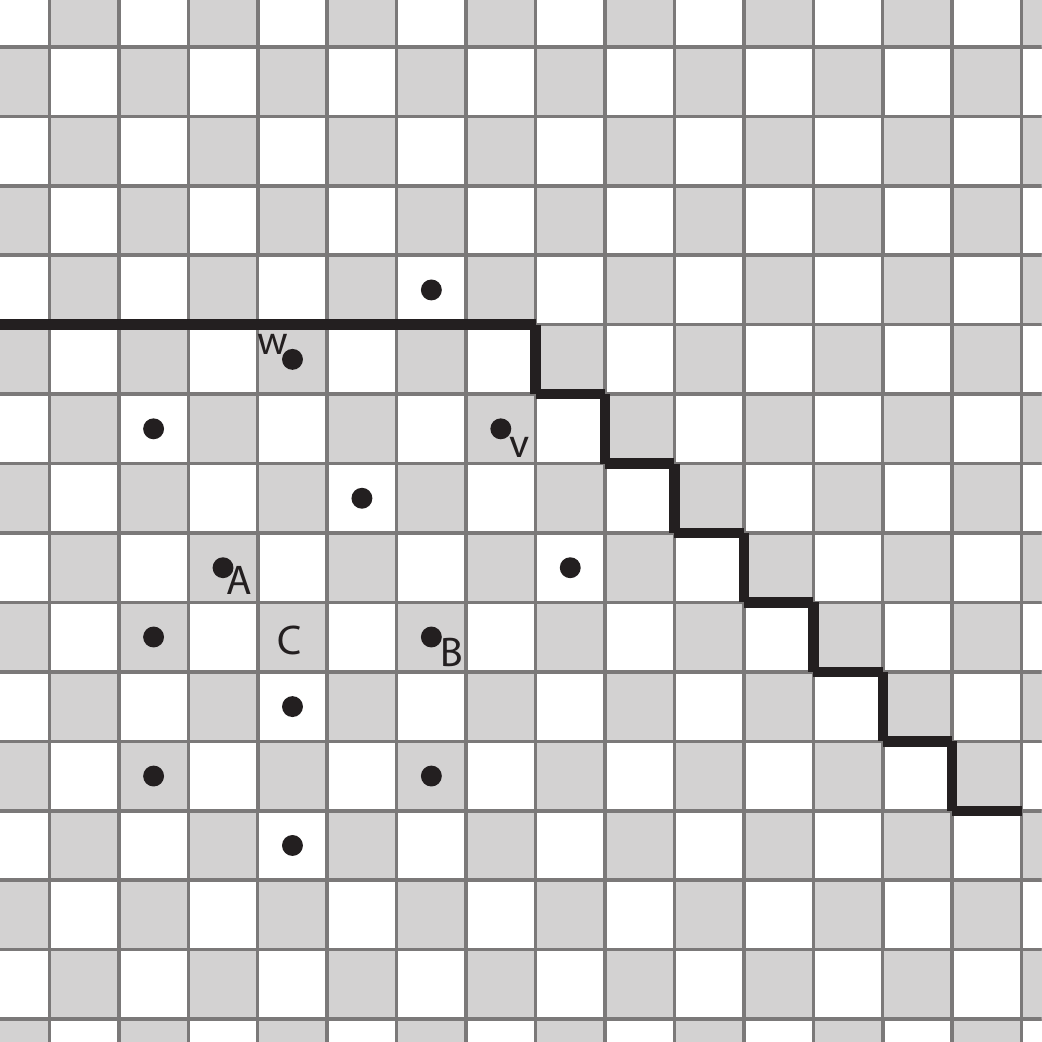}}
\end{center}
\caption{Case 4 (a) (ii).}
\label{picfourteen}
\end{figure}

(ii) $w$ is at $L$ of Figure \ref{spread1}. 
In this case, refer to Figure \ref{picfourteen}.
We can have $w$ connected to $A$ and $v$ connected
to $B$, both via enhanced diamond sites, with 
the enhancement at $C'$ switched off.

\begin{figure}%[htbp]
\begin{center}
\scalebox{1.2}{\includegraphics{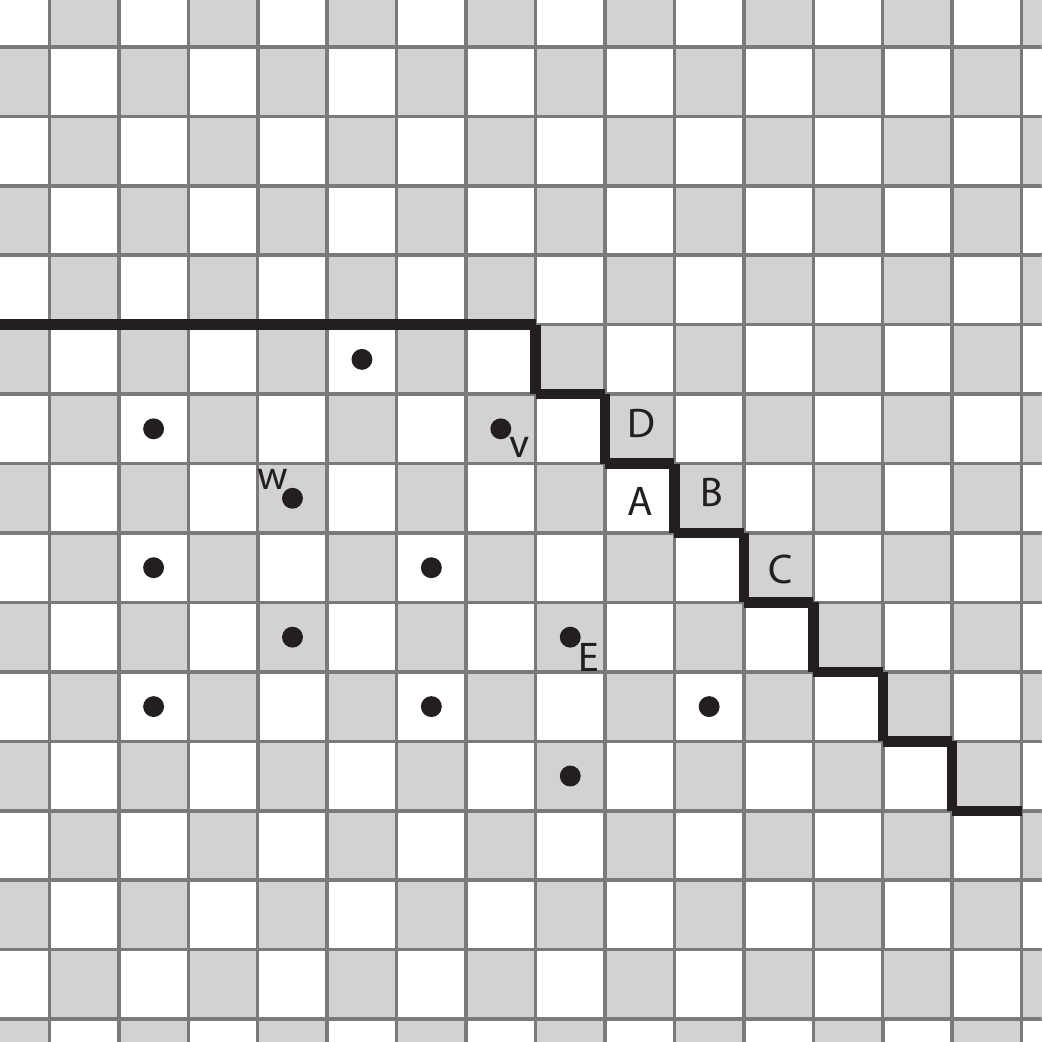}}
\end{center}
\caption{Case 4 (a) (iii).}
\label{spread3}
\end{figure}

 (iii) $w$ is at $K$ of Figure \ref{spread1}.
Then refer to Figure \ref{spread3}.
 We aim to have an occupied site at $E$ connected to $v$.
 This is fine as long as there is no site of $W$ at $B$ or $C$. If there is one at $C$ but not $B$ then we need to have an occupied odd site at $A$ and switch off the enhancement which we can do as we know there is no occupied site at $D$ as it would be joined to $v$ and $W$. If there is a site of $W$ at $B$ then it is not actually possible to have $y$ being $r$-pivotal as
 there is no way to get a path from $V$ into $D_r$ without joining up with $W$.
This is because $v$ is blocked from having a path
further into $G_r$, and there cannot be any other
point in $G_r$ connected to $V$ elsewhere, because
the paths in $W$  from locations in $G_r$ on
both sides of $v$ cut $v$ off 
from being path-connected to any other part of the boundary of $G_r$.

(b) $v$ is at $I$ of Figure \ref{spread1}.
 If $w$ is $3$ or more diagonals away then there is no problem.
 If $w$ is $4$ or more columns away then there is no problem. This just
 leaves two possibilities.

\begin{figure}%[htbp]
\begin{center}
\scalebox{1.2}{\includegraphics{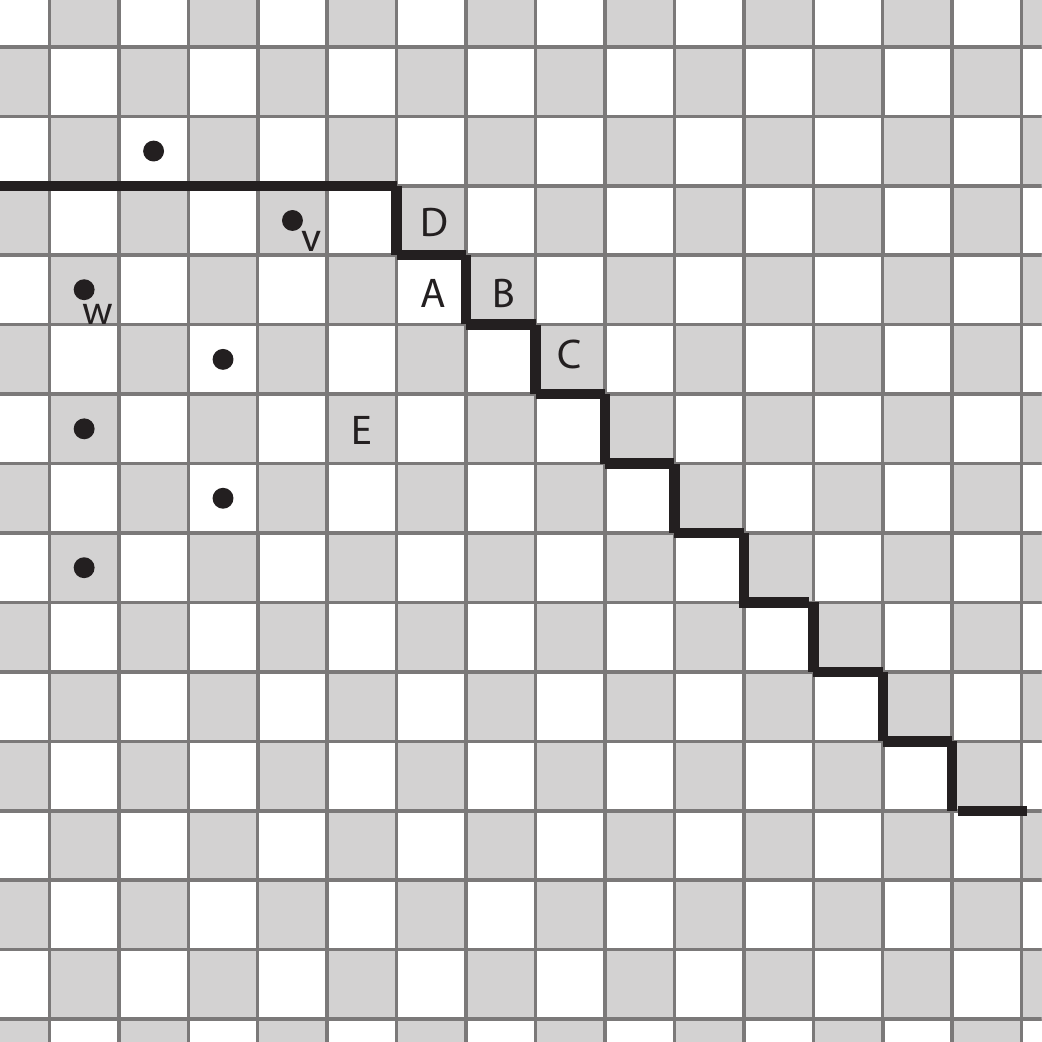}}
\end{center}
\caption{Case 4 (b) (i).}
\label{spread4}
\end{figure}

(i) If $w$ is at $O$ of Figure \ref{spread1},
then (see Figure \ref {spread4}) 
 this is akin to case (a) (iii) but just translated.

\begin{figure}%[htbp]
\begin{center}
\scalebox{1.2}{\includegraphics{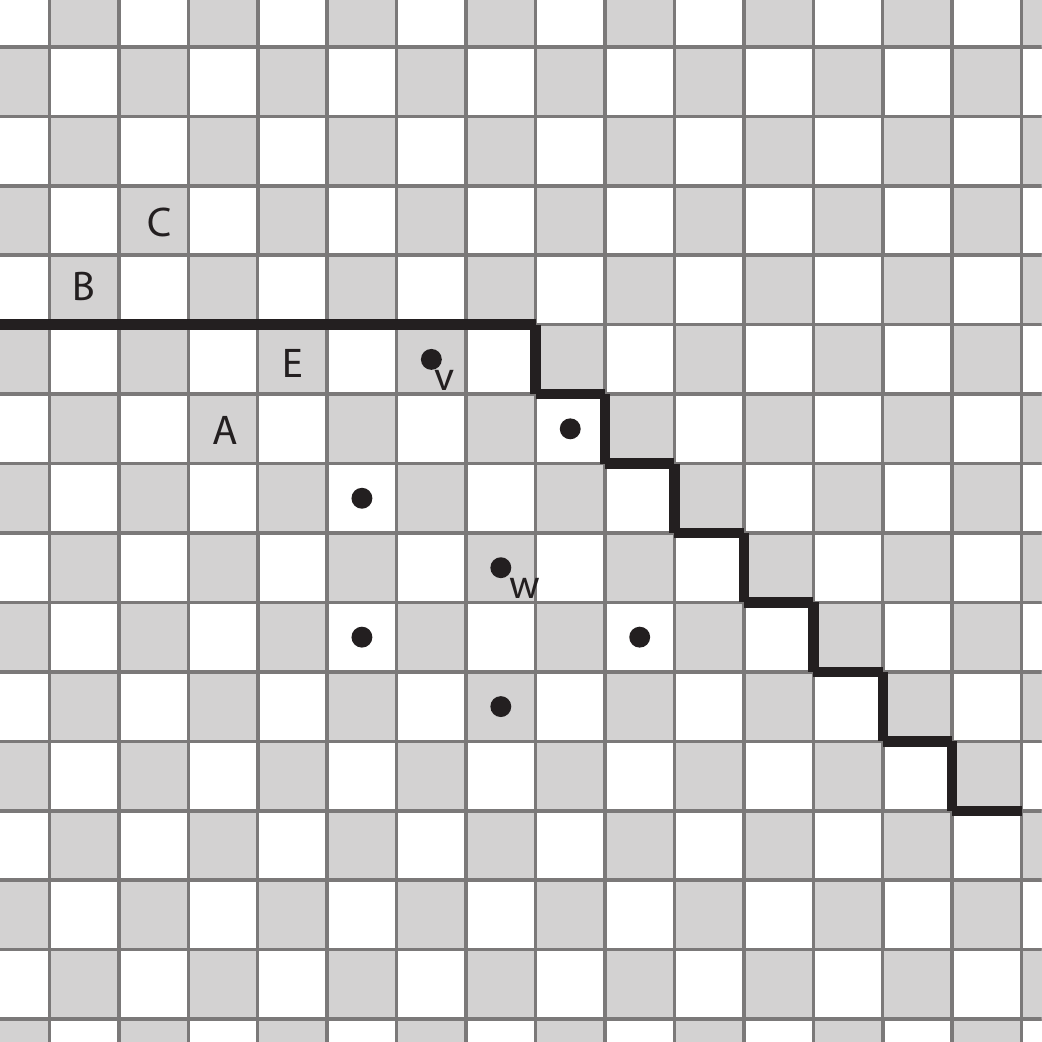}}
\end{center}
\caption{Case 4 (b) (ii).}
\label{spread5}
\end{figure}

(ii) If $w$ is at $N$
of Figure \ref{spread1},
then (see Figure \ref {spread5}) 
 we aim to have an occupied even site at $A$. We can do this unless there is an occupied site at $B$ which is in $W$. If this happens then we aim for an occupied even site at $E$ instead. This works so long as there is no occupied site at $C$ in $W$. So there is no problem unless there are occupied sites at both $B$ and $C$ in $W$. If this happens then it is not actually possible to have $y$ being $r$-pivotal as there is no way to get a path from $V$ into $D_r$ without joining up with $W$. 

Now consider the cases where $C_{2r}$ (and hence $G_r$)
 is not contained in $B(2n)$.
 First we look at the case where $C_{2r}$ intersects just the top
 edge of $B(2n)$. We define an octagonal region
$F_r$ as follows.
Start with the rectangular region $C_{2r} \cap B(2n)$.
Then remove triangular  regions at  the corners to 
make an octagon. The triangular regions are
of height 9 or 10, chosen in such a way
that the inner boundary consists of odd sites.
 We then argue as before using $F_r$ instead of $G_r$.
We have the sets $V$ and $W$ as before and the sites $v$ and $w$.
 If $v$ and $w$ are both well away from the edge of $B(2n)$ then we just have one of the cases we have already looked at. So we just consider the 
case where $v$ say is near the edge of $B(2n)$.
 However as it is on a diagonal of $F_r$ we can treat it as before and the path we create will stay inside $B(2n)$.

Now consider the case where $C_{2r}$ intersects the right hand edge
 of $B(2n)$, and define an octagonal  region $F_r $ inside
$B(2n)$ analogously to the previous case.
 In this case we just look at the set $V$ and site
 $v$ inside $F_{r}$ that is connected to the left of $B(2n)$.
 Inside $F_r$ we can then form a path from $v$ towards $y$ and a
 disjoint path from the right hand edge of $B(2n)$ towards $y$
 and ensure that $\ty'$ is $2$-pivotal. 

Finally we consider the case with $r > n/5$.
In this case, we can make a path of even sites  in from each boundary
of $B(2n)$ to $y$, together with a path of odd sites around the edge
of each of these paths and around the boundary of $B_n$. 
\hfill {$\Box$} \\

{\bf Proof of Proposition \ref{propP1P2}.} 
Assume $(\lambda,p) \in [0.5,1.5] \times [0.2,0.8]$.
Suppose 
 $E(r) $  occurs for the $S_x$ process.
Let the sets $Q_1, Q_2, Q_3$ be as in Lemma \ref{Qprop}.
Suppose also that
 $ E_1^{(2r+6,2r+10)} $
 occurs, and  we have $T_x \leq 1$ on all occupied
sites (for the $S_x$-process) in $C_{2r+4} \setminus G_r$
and $T_x >1$ on all blocked
sites (for the $S_x$-process) in $C_{2r+4} \setminus G_r$
(this is consistent with occurrence of event
 $ E_1^{(2r+6,2r+10)} $.)
Suppose also that $T_x \leq 1 $ for all 
the sites in
 $Q_1$ and $T_x > 1$ on all
 the sites in $ \Z^2 $ lying adjacent to $Q_1$,
and $T_{x'} < p$ for $x' \in Q_2$ and $T_{x'} > p$ for $x' \in Q_3$.
 Then using Lemma \ref{Rlem} we have that
  $y$ is $2$-pivotal for the $U_x$ process.
 This all occurs with probability at least 
$K_4^{-r}$ (given $E(r)$), for some finite positive constant $K_4$.
Therefore for all $y \in \Z^2 \cap B(n)$ and all $r \geq 20$ we have that
\bea
P_{2}(n,\lambda,\ty) \geq K_4^{-r} P[E(r)].
\label{1014c}
\eea
Hence by (\ref{1014b}) and (\ref{1111a}),
\bea
P_1(n,\lambda,y) =
 \sum_{r=0}^\infty P_{1,r,n}(y)
\leq 20 P[E(20)] + \sum_{r=20}^\infty 
\frac{K_2^r P[E(r)]}{\lfloor r/2 \rfloor !}
\nonumber \\
 \leq 20 K_4^{20} P_2(n,\lambda,\ty)
+ \sum_{r=20}^\infty \frac{(K_2K_4)^{r} P_2(n,\lambda,\ty)}{
\lfloor r/2 \rfloor !} = K_1 P_2(n,\lambda,\ty),
\eea
where $K_1$ is a finite constant independent of $\lambda$ and $p$,
as required.
%given by $K_1$ 
\hfill{$\Box$}

%\section{Comparison of 1-pivotal probabilities}
\section{Proof of Theorem \ref{thm1}}
\allco
%\eqco
In the preceding section we found a lower bound for
$P_1(n,\lambda,y)$ in terms of $P_2(n,\lambda,\ty)$, for
$y$ inside $B(2n)$. We now find a lower bound
for $P_1(n,\lambda,y)$ in terms of $P_1(n,\lambda,z)$
for $y$ outside $B(2n)$ and $z$ inside $B(2n)$.
Once we have this, we shall be able to quickly
complete the proof of Theorem \ref{thm1}.

We introduce more notation.
Let $\partial B(2n)$ be the set of even sites on the inner 
boundary of $B(2n)$.
For $x \in \Z^2 \setminus B(2n)$, let $z(x)$ be
the nearest site
 in $ \partial B(2n) $ to $x$
(here using  graph distance in $\Z^2$
 as our measure of distance). 
If there is a choice of two we take $z(x)$ to
 be the one clockwise from the other.
For $z \in B(2n)$, set $L_z := \{x \in \Z^2 \setminus B(2n): z(x)=z\}$. 
%We split the outside of $B(2n)$ into sets $L_z$ where a site $x$ is
% in $L_z$ if $z$ is 
\begin{prop}
 There exists a constant $K_5$ such that 
 for any $(\lambda,p) \in [0.5,1.5]\times [0.2,0.8]$ and any
$z \in \partial B(2n)$ and even $y \in L_z$ we have that
\bea
P_{1,r,n}(y) \leq \frac{
 K_5^r
%P_{1,r,n}(z)
P_{1}(n,\lambda,z)
 I_r(y)}{\lfloor r/2 \rfloor !},
~~~~ r \geq 20;
\label{1111b}
\\
P_{1,r,n}(y) \leq 
%P_{1,20,n}(z)
 K_5
P_{1}(n,\lambda,z)
 I_r(y),
~~~~ r \leq 20.
\label{1111c}
\eea
where $I_r(y) = 1$ if $y$ is within $r$ steps of $B(2n)$ and $I_r(y)=0$
 otherwise.
\lbl{lem1014}
\end{prop}

{\bf Proof.}
Assume $y$ is within $r$ steps of $B(2n)$; otherwise
 it cannot possibly be $(1,r)$-pivotal. 
The proof is very similar to that of Proposition \ref{propP1P2}.
We couple processes as before.
That is, we start with independent $S_x$
and $T_x$ arrivals processes, and
 define $U_x = U_x^{(2r+6,2r+10)} $ by (\ref{1014d}) as before.
%Define regions  $M := M,M_1,M_2, N, N_1,N_2$ and $R$ as before.

Given $y$, and given $r \in \N$,
define event $E(r)$ as in Section \ref{secpivcomp}.
%If $y$ is $(1,r)$-pivotal for the $S_x$ process then 
% event $E(r)$ must occur in the $S_x$ process. 
Although now $y$ lies outside $B(2n)$,
Lemma \ref{newElem} remains valid.
% We  assume that $E(r)$ occurs in the $S_x$ process. 

Let event $E_1 := E_1^{(2r+6,2r+10)}$  be
defined by (\ref{E1def})
as before.
By  Lemma \ref{Rlem},
%whatever the arrival 
%times are on $\Z^2 \setminus C_{2r+6}$, 
the state
 of all sites in $\Z^2 \setminus
C_{2r+6}$ will be
 the same in the $U_x$ process as in the $S_x$ process if $E_1 $ occurs.

%Suppose $y$ and $z$ are as in the statement of Proposition
%\ref{lem1014}, with $I_r(y)=1$.
Define the region $F_r$ as we did in the proof 
of Lemma \ref{Qprop} when there were boundary effects.
That is, take the intersection of $C_{2r} \cap B(2n)$ 
and smooth the corners to get an octagonal region with inner
 diagonal of length
9 or 10 consisting of odd sites. But if $C_{2r} \cap B(2n)$ shares
a corner with $B(2n)$, then do not smooth that particular
corner. Then $z$ will lie in the  region $F_r$.

Using Lemma \ref{newElem} and Lemma \ref{2Qprop} below,
 which is analogous to Lemma \ref{Qprop}, 
as in the proof of Proposition \ref{propP1P2}, 
we can find a constant $K_6$ such that
for $r \geq 20$ we have
$$
P_{1,r,n} (y) \leq  \frac{K_2^r P[E(r)]}{\lfloor r/2 \rfloor!}
 \leq  \frac{K_2^r K_6^rP_{1}(n,\lambda,z)}{\lfloor r/2 \rfloor!}
$$
which demonstrates (\ref{1111b}). 
In the case with $r \leq 20$, 
we  use Lemmas \ref{newElem} and \ref{2Qprop} 
to obtain
$$
P_{1,r,n}(y) \leq P[E(r)] \leq P[E(20)] \leq K_6^{20}
P_{1}(n, \lambda,z),
$$
yielding
 (\ref{1111c}).
%\ref{lem1014} in a similar manner to that of 
\hfill{$\Box$} \\

\begin{lemm}
There exists a constant $\alpha \in (0,\infty)$
with the following property.
Let $y,z$ be as above and assume $r \geq 20$.
If the event $E(r)$ occurs in the $S_x$ process,
then there exists
a stable  set $Q_1 \subset F_r \cap \Z^2$ having
no element  adjacent to the occupied sites of the $S_x$
process outside
$F_r$, and
disjoint sets $Q_2, Q_3$ of diamond sites 
inside $F_{r}$,
such that (i) each of $Q_1, Q_2, Q_3$ has at most $\alpha r$
elements, and (ii) if, 
in the $U_x$ process, all the sites in $Q_1$ are occupied, 
   all diamonds in $Q_2$ are black,  all the diamonds in $Q_3$ are white,
and all sites in $C_{2r+6} \setminus Q_r$ are in the
same state as for the $S_x$ process,
 then $z$ is $1$-pivotal
for the $U_x$  process.
\label{2Qprop}
\end{lemm}
{\bf Proof.}
Suppose $C_{2r}$ does not meet the left or right boundary
of $B(2n)$. If
 $E(r)$ occurs  there must be disjoint paths in the $S_x$ process up to 
$\Z^2 \cap C_{r+7}$ within $B(2n)$ from each side of $B(2n)$. 
By similar arguments to those in the proof of Lemma
\ref{Qprop}, 
  we can obtain the event that $z$ is $1$-pivotal 
for the $U_x$ process, by specifying $O(r)$ vertices
%of the $T_x $ process
 to be occupied.

Suppose  $C_{2r}$ meets the right  boundary of $B(2n)$.
 Then if $E(r)$ occurs  there must be a path in the $S_x$ process up to 
$C_{r+7}$ within $B(2n)$ from the left side of $B(2n)$. 
Hence there is such  path from the left boundary of
$B(2n)$ to the boundary of $F_r$.
By similar arguments to before, 
%if also $E_1 \cap E_2$ occurs then
 we can obtain the event that $z$ is $1$-pivotal 
for the $U_x$ process, by specifying $O(r)$ vertices
%of the $T_x $ process 
to be occupied so as to extend the existing
path to $z$, and creating a disjoint path from the right hand
edge of $B(2n)$ to $z$.

The case where $C_{2r}$ meets the left  boundary of $B(2n)$ is
treated analogously.
\hfill{$\Box$} \\

{\bf Proof of Proposition \ref{keyprop}.}
By Proposition  \ref{lem1014}, there are constants $K_7$, $K_8$
 such  that
for any $z \in \partial B(2n)$,
\bean
\sum_{y \in L_z: y ~{\rm even}} P_1(n,\lambda ,y)
= \sum_{y \in L_z} \sum_{r=0}^\infty
P_{1,r,n}(y)
\\
\leq \sum_{r=0}^{19} K_7 P_{1}(n,\lambda,z) +
 \sum_{r=20}^\infty \frac{P_{1}(n,\lambda,z)
K_5^r K_7 r^2}{\lfloor r/2 \rfloor !} \leq
K_8 P_1(n,\lambda,z).  
%(**)
\eean
%and therefore there exists $\delta_3$ independent of 
%$n$ and $x$ such that for all $z \in \partial B(2n)$,
%\bean
%P_1(z) \geq \delta_3 \sum_{x \in L_z} P_1(x).
%\eean
Summing over  $z \in \partial B(2n)$, we obtain that
%have that
%there exists $\delta_2(\lambda,p)$ which is independent of $n$
% and strictly positive and continuous on $[0.5,1.5] \times (0,1)$ such that
\[
\sum_{y \in \Z^2 \setminus B(2n): y ~{\rm even} } P_1(n,\lambda,y)
\leq 
K_8
\sum_{z \in B(2n) \cap \Z^2: z ~ {\rm even}} P_1(n,\lambda,z). 
\]
Putting this   together with Proposition \ref{propP1P2}
gives for some $K_9$ that
\[
\sum_{y \in \Z^2: y~ {\rm even}} P_1(n,\lambda,y)
\leq  K_9
\sum_{z \in \Z^2:z' \in B(2n)} P_2(n,\lambda,z). 
\]
Hence by Proposition \ref{Russoprop},
\[
\frac{\partial h(n,\lambda,p)}{\partial \lambda}
\leq 
K_9 \frac{\partial h(n,\lambda,p)}{\partial p} , ~~~~~
(\lambda,p) \in [0.5,1.5] \times [0.2,0.8].
\]
We also know from (\ref{1024a}) 
that $h(n,1,0.5) = 0.5$, so looking at
 a small box around $(1,0.5)$ we can find
 $\epsilon >0$ such that for all $n$,
we have
 $h(n,1-\epsilon,1) \geq h(n,1,0.5) = 0.5$.
 Therefore taking $\mu = 1-\epsilon$ we have satisfied (\ref{eqstar}). 
\hfill{$\Box$} \\ 

With Proposition \ref{keyprop} proven,
our proof of Theorem \ref{thm1} is now complete by
the arguments already given in Sections 
%We already proved $\lambda_c< 10$ in Section
 \ref{secintro} and \ref{secdual}. \\ 
%We can then obtain the lower bound $\lambda_c >1$  by the argument
%at the end of Section \ref{secdual}, and 

{\bf Acknowledgements.}
We thank Martin Zerner and G\"unter Last for
some helpful discussions pertaining to this
paper.

\end{document}